\documentclass{article}

\usepackage{amsmath}
\usepackage[francais]{babel}
\usepackage{amssymb}
\usepackage[applemac]{inputenc}
\usepackage{mathrsfs}

\def\noi {\noindent}
\def\Q{\mathop{\bf Q}}
\def\F{\mathop{\bf F}}

\def\Z{\mathop{\bf Z}}

\def\R{\mathop{\bf R}}
\def\C{\mathop{\bf C}}
\def\s{\mathop{\bf S}}
\def\H{\mathop{\bf H}}

\def\deg{\mathop{\rm deg}}
\def\ab{\mathop{\rm ab}}
\def\dim{\mathop{\rm dim}}
\def\Sym{\mathop{\rm Sym}}
\def\Alt{\mathop{\rm Alt}}
\def\det{\mathop{\rm det}}
\def\rg{\mathop{\rm rg}}
\def\rgr{\mathop{\rm rgr}}
\def\GL{\mathop{\rm GL}}

\def\SO{\mathop{\rm SO}}
\def\O{\mathop{\rm O}}
\def\Sp{\mathop{\rm Sp}}

\def\SU{\mathop{\rm SU}}

\def\Aut{\mathop{\rm Aut}}

\def\Inv{\mathop{\rm Inv}}

\def\cc{\mathscr{C}}

\def\bb{\mathscr{B}}

\def\hh{\mathscr{H}} 

\def\sm{\raisebox{1.5pt}{~\rule{5pt}{1.3pt}~}}

\begin{document}
 
  \begin{center} 
  
  Groupes de Coxeter finis: involutions et cubes 
 
 \bigskip 
   Jean-Pierre Serre
   
   \end{center}
   
   \vspace{10mm} 
   
     Le texte qui suit passe en revue les propriétés des groupes de Coxeter finis
     qui sont les plus utiles pour la compréhension de leurs {\it invariants cohomologiques}, au
  sens de [Se 03]. La plupart de ces propriétés se trouvent déjà
     dans la littérature (par exemple [Bo 68], [Ca 72], [Sp 74], [Sp 82], [Ka 01], [DPR 13]),
     mais il a paru commode de les rassembler de façon systématique, et de les compléter sur quelques points.
     
       Il s'agit surtout des {\it involutions} ($\S$2), car ce sont leurs classes de conjugaison qui paramètrent de façon naturelle les invariants cohomologiques
 du groupe, cf. [Se 18]. Elles interviennent le plus souvent par l'intermédiaire des sous-groupes appelés ``{\it cubes}'': sous-groupes abéliens engendrés par des réflexions ($\S$4). Ces groupes jouent un rôle analogue à celui des groupes de Sylow, cf. par exemple th.4.16. Leur intérêt
 pour la détermination des invariants cohomologiques provient du ``{\it splitting principle}'': sous certaines conditions techniques, un invariant cohomologique est nul si ses restrictions à tous les cubes sont nulles, cf. [Se 03], [Hi 10], [Se 18],
 [Hi 20], [GH 21].
 
   Le $\S$5 décrit les différents types de groupes irréductibles,
   en insistant sur les propriétés de leurs involutions et de leurs cubes, notamment pour les types $B_n$, $D_n$, $E_7$
   et $E_8$. Il contient aussi une liste d'inclusions entre différents types
   qui se révèle utile pour prouver certains cas du ``splitting principle'',
   comme nous le montrerons ailleurs. Le $\S$6 est consacré à la construction, et aux propriétés, de certains groupes de Coxeter de rang 4, notamment ceux
   de type $F_4$ et $H_4$, cf. th.6.12.

  \bigskip

 {\bf $\S$1. Groupes de Coxeter finis}. 
 
 \medskip
 
\noi  1.1. Dans ce qui suit, $V$ désigne un espace vectoriel réel de dimension finie,
 et $G$ un sous-groupe fini de $\GL(V)$. On suppose que $G$ est engendré par des réflexions, c'est-à-dire par des éléments d'ordre 2 fixant un hyperplan (rappelons que ``fixer'' un sous-ensemble signifie fixer tout point de ce sous-ensemble). C'est donc un groupe de Coxeter fini, et tout groupe de Coxeter fini s'obtient
 de cette façon, cf. [Bo 68], $\S$V.4. Un tel couple $(V,G)$ sera appelé un {\it couple de Coxeter}.
 
   On munit $V$ d'un {\it produit scalaire} $\langle x,y\rangle$ invariant par $G$
   tel que la forme quadratique $x \mapsto \langle x,x\rangle$ soit définie positive.
   La notion d'orthogonalité sera relative à ce produit scalaire.

 \medskip
\noi 1.2.   Soit $V_0$ le sous-espace de $V$ fixé par $G$, et soit $V_1$ le sous-espace de $V$ engendré par les $(g-1)x, g \in G, x\in V$. Alors $V = V_1 \oplus V_0$; de plus $V_1$ est l'orthogonal de $V_0$.
   
   On dira que le couple $(V,G)$ est {\it réduit} si $V_0 = 0$. Nous ne considérerons la plupart du temps que des couples réduits, mais il y a certains cas (tel celui du groupe symétrique $\Sym_n$) où l'emploi d'un $V$ non réduit (à savoir $\R^n$) est plus commode.
   
   La dimension de $V_1$ est appelée le {\it rang} de $(V,G)$, et notée $\rg(G)$.
   
    \medskip
 \noi 1.3.  On dit que $G$ est {\it cristallographique} s'il stabilise un réseau de $V$,
 ou, ce qui revient au même, s'il stabilise une $\Q$-structure sur $V$.  D'après
 [Bo 68], $\S$VI.2, prop.9, cela équivaut à dire que $G$ est le {\it groupe de Weyl} d'un système de racines de $V$, que l'on peut supposer ``réduit'' au sens de [Bo 68],
  $\S$VI.1.4.    
    
 \medskip
 \noi 1.4. Un couple $(V,G)$ est dit {\it irréductible} si $G \ne 1$ et si $V$ est un $G$-module irréductible ; cela entraîne que $(V,G)$ est réduit. On dit aussi que $G$ est irréductible.
 
  Tout couple réduit$(V,G)$ se décompose de façon unique en $V = \oplus_i (V_i,G_i)$ 
  avec $G = \prod_i G_i$, les $(V_i,G_i)$ étant irréductibles cf. [Bo 68], $\S$V.3.7. On dit alors que les $G_i$ sont les facteurs irréductibles de $G$. Les $V_i$
  sont orthogonaux entre eux.
    
  On trouvera dans [Bo 68], $\S$VI.4, th.1, la liste des $G$ irréductibles.
  Il y a les cristallographiques:

  \smallskip

  $A_n \ (n \geqslant 1),\ B_n \ (n \geqslant 2), \ D_n \  (n\geqslant 3), \ E_6, \ E_7, \ E_8, \ F_4, \ G_2$,
  
  \smallskip
  \noi et aussi: $H_3, \ H_4, \ I_2(m)$ avec $ m \geqslant 3$ (groupe diédral d'ordre $2m$). 
  
  \smallskip
    Dans chaque cas, l'indice dénote le rang. Les seuls isomorphismes entre 
  éléments de cette liste sont $D_3 \simeq A_3$, $A_2 \simeq I_2(3), B_2 \simeq I_2(4)$ et $G_2 \simeq I_2(6).$ Noter l'absence de la série $C_n$ car les groupes de Weyl correspondants sont isomorphes à ceux de type $B_n$.

    Noter aussi que le groupe diédral $I_2(2)$ d'ordre 4 est un groupe de Coxeter
    réductible, de type $A_1 \times A_1$.
    
   \medskip
  \noi 1.5. Un sous-groupe $H$ de $G$ sera appelé un {\it $\cc$-sous-groupe} ({\it reflection subgroup} en anglais) s'il est engendré par des réflexions de $G$, autrement dit si le couple $(V,H)$ est un couple de Coxeter.

 \smallskip
  \noi {\bf Théorème 1.6.} {\it Si $X$ est une partie de $V$, le sous-groupe $G_X$
  de $G$ formé des éléments fixant $X$ est un $\cc$-sous-groupe de $G$}.
  
  \smallskip
  \noi {\it Démonstration.} C'est un cas particulier de [Bo 68], $\S$V.3, prop.2. 
  
  \smallskip
  \noi {\it Remarque}. Les sous-groupes de $G$ de la forme $G_X$ sont appelés
  les {\it sous-groupes paraboliques} de $G$. Soit $S$ une base de Coxeter de $G$ (i.e. un ensemble de réflexions tel que le couple $(G,S)$ soit un système de Coxeter
  au sens de [Bo 68], $\S$IV.1.3). Alors un sous-groupe $H$ de $G$ est parabolique si et seulement si l'un de ses conjugués est engendré par une partie de $S$; cela résulte par exemple de [Bo 68], $\S$V.3.3. 
    
  \medskip

  \noi 1.7. {\it Centralisateur d'une réflexion.}
  
  \smallskip
    Soit $s$ une réflexion et soit $D_s$ la droite de $V$ formée des
    éléments $x$ tels que $sx=-x.$ Soit $G_s^{^+}$ le fixateur de $D_s$, et soit
    $G_s$ le centralisateur de $s$ dans $G$. 
    
    \smallskip
    \noi {\bf Proposition 1.8}.  \  $G_s = \{1,s\} \times G_s^{^+}.$ 
    
    \smallskip
      \noi   {\it Démonstration.} L'inclusion $G_s \supset \{1,s\} \times G_s^{^+}$ est claire. Si $g \in G_s$, on a $gD_s = D_s$. Or un automorphisme d'ordre fini d'une droite réelle est égal à $\pm 1$. On a donc, soit $g \in G_s^{^+}$, soit $sg \in G_s^{^+}$. D'où la proposition.

\smallskip  

\noi {\bf Corollaire 1.9}. {\it Les groupes  $G_s^{^+}$ et  $G_s$ sont des $\cc$-sous-groupes de $G$. }

\noi Cela résulte du th. 1.6.
      
      \medskip
      
  \noi 1.10. Puisque $G$ est un groupe de Coxeter, il peut se définir par générateurs
  et relations comme expliqué dans [Bo 68], $ \S$IV.1. L'ensemble générateur $S$ se compose de $n$ réflexions $s_1,...,s_n$, où $n = \rg G$. Les relations sont de la forme $(s_is_j)^{m_{ij}}=1$, où les $m_{ij}$ sont des entiers $\geqslant 1$, égaux à 1 si et seulement si $i=j$, et tels que $m_{ij}=m_{ji}$ pour tout $i,j$. Soit $M_G$ l'ensemble des $m_{ij} > 2$. 
  
  Si $G = \prod G_\lambda$ est la décomposition de $G$ en facteurs irréductibles, on a $M_G = \bigcup M_{G_\lambda}$. On a $M_G = \varnothing$ si et seulement si tous les facteurs irréductibles de $G$ sont de type $A_1$, autrement dit si $G$ est un {\it cube}, au sens du $\S$4.
  
  \smallskip
  Voici la liste des $M_G$ pour $G$ irréductible:
   
   \smallskip
  
  \noi    $M_G =\{3\}$ \quad \ si $G$ est de l'un des types $A_n$ ($n\geqslant 2)$, $D_n$, $E_6$, $E_7$, $E_8$; 
  
    \noi $M_G = \{3,4\}$ \ si $G$ est de type $B_n$ ou $F_4$;
 
 \noi $M_G = \{3,6\}$ \ si $G$ est de type $G_2$; 

 \noi $M_G=\{3,5\}$ \ si $G$ est de type $H_3$ ou $H_4$;

 \noi $M_G=\{m\}$  \quad si $G$ est de type $I_2(m)$ .

\smallskip

Le cas $M_G = \varnothing$ ou $\{3\}$ est celui où le graphe de Coxeter ne contient pas de traits multiples (il est ``{\it simply laced} '').

\smallskip

\noi {\bf Théorème 1.11}. {\it Soit $n$ un entier $\geqslant 3$.  Pour qu'il existe un $\cc$-sous-groupe de $G$ isomorphe à $I_2(n)$, il faut et il suffit que $n$ divise l'un des éléments de $M_G$.}

  \noi  (La condition ``contenir un $\cc$-sous-groupe isomorphe à $I_2(n)$'' équivaut à dire qu'il existe deux réflexions dont le produit est d'ordre $n$.)

\smallskip

\newpage

\noi {\it Démonstration} (suggérée par le rapporteur). Si $n$ divise un élément $m$ de $M_G$, il existe des réflexions $s_1,s_2$ de $G$ telles que $s_1s_2$ soit d'ordre $m$; dans ce cas, le produit des réflexions $s_1$ et $s_2(s_1s_2)^{m/n}$ est d'ordre $n$, ce qui montre que $G$ contient un $\cc$-sous-groupe isomorphe à $I_2(n)$. Inversement, soit $H$ un sous-groupe de $G$ engendré par deux réflexions dont le produit est d'ordre $n$. Soit $X$ le sous-espace de $V$ fixé par $H$, et soit $G_X$
le fixateur de $X$. Le groupe $G_X$ contient $H$, et il opère fidèlement sur l'orthogonal $Y$ de $X$
dans $V$, qui est de dimension $2$; c'est donc un sous-groupe parabolique de rang 2 de $G$; il
est diédral d'ordre $2m$, où $m$ est un multiple de n. Puisqu'il est parabolique, l'un de ses conjugués
est engendré par un couple $s_i,s_j$ d'éléments de la base $S$ de 1.10. On en conclut que $m$ appartient à $M_G$, donc que $n$ divise un élément de $M_G$.

%\noi {\it Démonstration}. Il suffit de le prouver lorsque $G$ est irréductible. Le cas où 
%$G$ est cristallographique résulte de la description des diverses possibilités des couples de racines donnée dans [Bo 68], $\S$VI.1.3. Le cas du type
%$I_2(m)$ résulte de ce que $I_2(m)$ contient $I_2(n)$ si $n|m$. Le cas du type
%$H_3$ se vérifie sur la description de ce groupe comme $\Alt_5 \times \{\pm1\}$.
%Le cas de $H_4$ résulte du précédent puisque le centralisateur d'une réflexion
%est de type $A_1 \times H_3$.

\smallskip

\noi {\bf Corollaire 1.12.} {\it On a $M(H) \subset M(G)$ pour tout  $\cc$-sous-groupe $H$ de $G$.}

\medskip

\noi 1.13. {\it Groupes de type impair.} 

\smallskip
  
 Nous dirons que le type de $G$ est {\it impair} si tous les éléments de $M_G$ sont impairs, autrement dit si $G$ n'a aucun $\cc$-sous-groupe de type
  $I_2(m)$ avec $m$ pair $\geqslant 4$. 
  
  Les types irréductibles impairs sont : $A_n, D_n, E_6, E_7, E_8, H_3, H_4$ et
  les $I_2(m)$ pour $m$ impair.
  
  Un groupe de Coxeter est de type impair si et seulement si tous ses facteurs irréducibles le sont.
  
   \smallskip
  
  \noi {\bf Proposition 1.13}. {\it Si  $G$ est de type impair, il en est de même de ses $\cc$-sous-groupes.}
  
  \smallskip
  
  Cela résulte du cor.1.12.
  
  \smallskip
  
  \noi {\bf Proposition 1.14}. {\it Supposons $G$ irréductible. Si $-1 \not \in G$,
   alors $G$ est de type impair.}
   
   \smallskip
   
   En effet, $-1 \not \in G$ entraîne que $G$ est de l'un des types $I_2(m) \ (m$ impair), $A_n \ (n > 1)$, $D_n \ (n$ impair) ou $E_6$, et $G$ est donc de type impair.
   
   \smallskip
  
    \noi {\bf Théorème 1.16}. {\it Supposons $G$ irréductible. Les deux propriétés suivantes sont équivalentes}:
  
   (i) {\it $G$ est de type impair.}
   
   (ii) {\it Les réflexions de $G$ sont conjuguées entre elles.}
 
 \smallskip
 \noi  {\it Démonstration de} (i)$ \ \Rightarrow $ (ii). 
 
 Supposons que (i) soit vérifiée. Soient $S = \{s_1,...,s_n\}, m_{ij}$  comme dans 1.10. Par hypothèse, les $m_{ij}$ avec $i \neq j$ sont, soit égaux à 2, soit  impairs; dans le cas où $m_{ij}$ est impair, $s_i$ et $s_j$ sont conjugués dans le groupe diédral $\langle s_i,s_j\rangle$, donc dans $G$. Cela montre que les $s_i$ appartenant à la même composante connexe du graphe de Coxeter sont conjugués; or le graphe d'un groupe de Coxeter irréductible est connexe ([Bo 68], $\S$IV.1.9). Les réflexions $s_i$ sont donc conjuguées entre elles; comme toute réflexion de $G$ est conjuguée à l'une des $s_i$ ([Bo 68], $\S$V.3.2, cor. au th.1), cela entraîne (ii).

   \noi  {\it Démonstration de} (ii)  $ \Rightarrow$  (i). 

  Supposons (ii). Si $G$ n'était pas de type impair, il serait, soit de type $I_2(m)$, avec $m$ pair,
  soit de type $B_n$ ou $F_4$. Le cas de $I_2(m)$, $m$ pair, est impossible: ce groupe est engendré par deux réflexions qui ne sont pas conjuguées. Les autres cas 
  correspondent à des systèmes de racines où les racines peuvent avoir
  deux longueurs différentes (il y a des ``courtes'' et des ``longues''), donc définissent des réflexions non conjuguées.

  \bigskip
 
  {\bf $\S$2. Involutions}. 
 
 \smallskip
 
 \noi Soit $(V,G)$ un couple de Coxeter.
 
 \smallskip
 
 \noi 2.1.  Un élément $u$ de $G$ est appelé une {\it involution} si  $u^2=1$. Le {\it degré} d'un tel $u$ est la multiplicité de la valeur propre -1, ou, ce qui revient au même, la codimension du sous-espace $V^{^+}_u $ de $V$ fixé par $u$. La seule involution de degré $0$ est $u=1$. Les involutions de degré 1 sont les {\it réflexions}.
 
   Si $G = \prod G_i$ est la décomposition de $G$ en facteurs irréductibles, un
   élément $u = (u_i)$ de $G$ est une involution si et seulement si chaque $u_i$ est une involution\footnote{Ce ne serait pas le cas si l'on avait défini les involutions comme d'habitude en théorie des groupes, en demandant que leur ordre soit égal à $2$; il est commode d'accepter que $1$ soit une involution. Dans la correspondance
   ``classes d'involutions $\Longleftrightarrow$  invariants cohomologiques'' de [Se 18], l'involution $1$ correspond à un invariant constant.}. On a alors $\deg(u) = \sum \deg(u_i).$

 \smallskip
   \noi 2.2. Il est souvent commode de {\it ramener les résultats généraux sur les involutions au cas particulier de l'involution}  $ -1.$ Cela se fait de la manière suivante.
   
   Soit $u$ une involution, et soit $V = V^{^+}_u \oplus V^{^-}_u$ la décomposition de $V$ suivant les valeurs propres \ $+1,-1$ \ de $u$; on a $\deg(u) = \dim V^{^-}_u$.  Soit $G_u^{^-}$ le sous-groupe de $G$ formé des éléments qui fixent $V^{^+}_u$; d'après (1.6), c'est un $\cc$-sous-groupe de $G$, et $u$ est l'involution ``$ -1 $''  de $G_u^{^-}$, dans sa représentation naturelle sur~$V^{^-}_u$.
   
   \smallskip
   
   \noi 2.3. {\it Centralisateur d'une involution.}
   
   \smallskip
    Conservons les notations de 2.2. Soit $G_u^{^+}$ le fixateur de $V^{^-}_u$, et soit $G_u$ le centralisateur de $u$. On a
 $ G_u^{^-} \times G_u^{^+} \  \subset \ G_u. $ Les groupes $G_u^{^-}$ et $G_u^{^+}$ sont
 engendrés par des réflexions (cf. 1.6). De façon plus précise, si $s$ est une réflexion appartenant à $G_u$, on a $s \in G_u^{^+}$ si $\deg(us)=\deg(u)+1$ et
 $s \in G_u^{^-}$ si $\deg(us)=\deg(u)-1$. En particulier, $ G_u^{^-} \times G_u^{^+}$ {\it est le plus grand $\cc$-sous-groupe de $G_u$.} Lorsque $\deg(u) > 1$, il peut se faire que $G_u$ soit strictement plus grand que $ G_u^{^-} \times G_u^{^+}$, autrement dit que $G_u$ {\it ne soit pas un $\cc$-sous-groupe de $G$}.
   
   \smallskip
 \noi {\it Exemple.} Prenons $G$ de type $A_3$, donc isomorphe à $\Sym_4$; prenons pour $u$ une involution de rang $2$, autrement dit un produit de deux transpositions à supports
   disjoints. On a alors $G_u^{^+}= 1$ et $G_u^{^-}$ est de type $A_1 \times A_1$,
   alors que $G_u$ est un groupe diédral d'ordre $8$, donc n'est pas un $\cc$-sous-groupe de $G$. On verra au $\S$3.15 d'autres
   exemples du même genre.

 \smallskip
 \noi {\it Remarque}. Lorsque $u$ est une involution de rang 2, le groupe $G_u^{^-}$ est un groupe
   diédral $I_2(m)$ contenant $-1$, donc tel que $m$ soit pair (éventuellement égal à 2). L'entier $m$ est un invariant de la classe de conjugaison de $u$; d'après 5.11, $m/2$ est le nombre de décompositions de $u$ comme produit de deux réflexions.
   
   \smallskip
\noi {\bf Proposition 2.4}. {\it Si  $s_1,...,s_n$  sont des réflexions distinctes qui commutent deux à deux, le groupe qu'elles engendrent est d'ordre $2^n$ et leur produit $u = s_1\!\cdots\! s_n$ est une involution de degré  $n$}.

\smallskip
\noi {\it Démonstration}.  Soit $D_i$ la droite associée à $s_i$ comme dans 1.7.
Les $s_j$, $j \neq i$, fixent $D_i$, donc les $D_j$ sont contenus dans l'hyperplan fixé par $s_i$. Cela entraîne que $s_i$ n'est égal à aucun produit des $s_j$; l'ordre
du groupe $\langle s_1,...,s_n\rangle$ est donc $2^n$, et les droites $D_i$
sont orthogonales deux à deux. D'où le fait que le degré de $u$ est~$n$. 

 \smallskip
 
 \noi {\bf Proposition 2.5}. {\it Toute involution est produit de réflexions distinctes commutant deux à deux.}
 
  \smallskip
   
  \noi {\it Démonstration.} Soit $u$ une involution. On raisonne par récurrence sur  $\deg(u)$. On peut supposer que $\deg(u) \geqslant 2$, et aussi (d'après 2.2) que $u=-1$. Soit $s$ une réflexion de $G$; l'involution $u'=-s$ est de degré $\deg(u)-1$, et on lui applique l'hypothèse de récurrence.
  
  \medskip
  
  \noi 2.6. {\it Classes d'involutions}.
  Ce qui est important pour la détermination des invariants cohomologiques du groupe $G$, c'est l'ensemble $\Inv$ des {\it classes de conjugaison des
 involutions}. On trouve leur description (pour les groupes de type cristallographique)  
 dans [Ca 72], et [Ka 01], $\S\S$27-28. Nous reviendrons là-dessus au $\S$5.
 
 \medskip

\noi 2.7. {\it Nombre des classes d'involutions de degré donné}.

\smallskip

  Pour $n \geqslant 0$, soit $\Inv_n$ l'ensemble des classes de conjugaison des involutions de degré $n$. L'ensemble $\Inv$ est réunion disjointe des $\Inv_n$. Posons\footnote{Si $X$ est un ensemble fini, on note $|X|$ le nombre de ses
  éléments.}
  $h_n = |\Inv_n|$ et soit $h_G(t) = \sum h_nt^n$. On a $h_{G_1\times G_2}(t)= h_{G_1}(t)\cdot h_{G_2}(t).$ Il suffit donc de connaître les $h_G(t)$ lorsque $G$ est irréductible. Voici leurs valeurs:  
 
 \medskip
 
$\bullet$  $A_n: 1 + t + \cdots + t^{[(n+1)/2]}.$
  
 \smallskip
 
 \ \ Exemple - $A_3: 1 + t + t^2.$
 
 \medskip

$\bullet$  $B_n, \ n$ pair: $(1 + t +\cdots + t^{n/2})^2.$
 
 \smallskip

  \ \ Exemple - $B_{10}: 1+2t+3t^2+4t^3+5t^4+6t^5+5t^6+4t^7+3t^8+2t^9+t^{10}.$
    
  \medskip
  
$\bullet$   $B_n, \ n$ impair: $ (1+t+ \cdots + t^{(n-1)/2})(1+t+ \cdots + t^{(n+1)/2})$.
  
  \smallskip

\ \ Exemple - $B_9: 1+2t+3t^2+4t^3+5t^4+5t^5+4t^6+3t^7+2t^8+t^9.$
 
  \medskip
  
 $\bullet$  $D_n, \ n$ impair: $(1+t + \cdots + t^{(n-1)/2})(1+t + \cdots + t^{(n+1)/2})/(1+t)$.
  
  \smallskip
  
  \ \ Exemple - $D_{11}: 1+t+2t^2+2t^3+3t^4+3t^5+3t^6+2t^7+2t^8+t^9+t^{10}.$
  
   \medskip

 $\bullet$   $D_n, \ n$ pair: \ $t^{n/2} + (1+t + \cdots + t^{n/2})(1+t + \cdots + t^{1+n/2})/(1+t)$.

\smallskip
\ \ Exemple - $D_4: 1+t+3t^2+t^3+t^4.$

 \ \ Exemple - $D_{10}: 1+t+2t^2+2t^3+3t^4+4t^5+3t^6+2t^7+2t^8+ t^9+t^{10}$.
  
  \medskip
 
$\bullet$  $E_6: 1+t+t^2+t^3+t^4.$ 
 \medskip
 
 $\bullet$  $E_7:  1+t+t^2+2t^3+2t^4+t^5+t^6+t^7.$  
 \medskip

$\bullet$  $E_8: 1+t+t^2+t^3+2t^4+t^5+t^6+t^7+t^8$.
 
 \medskip
$\bullet$  $F_4: 1+2t+2t^2+2t^3+t^4.$
 
 \medskip
$\bullet$  $G_2: 1+ 2t +t^2.$
 
 \medskip
 
$\bullet$  $H_3: 1+t+t^2+t^3.$
 
 \medskip
 
 $\bullet$ $H_4: 1+t+t^2+t^3+t^4.$
 
  \medskip
 
  $\bullet$ $I_2(m): 1+t$ \ \ si $m$ est impair; $1+2t+t^2$ \ \ si $m$ est pair.
 
  \medskip
 
  Le degré $d$ de $h_G(t)$ est le {\it rang réduit} rgr($G$) de $G$, i.e. le degré des involutions maximales de $G$, 
 cf. 3.1.  On a $h_n = h_{d-n}$ pour tout $n$, cf. 3.11; autrement dit, $h_G$
 est un {\it polynôme réciproque}. On a aussi  $h_n \leqslant h_{n+1}$ pour $n < d/2$: cela se vérifie cas par cas lorsque $G$ est irréductible, et le cas général s'en déduit
 par un argument combinatoire; il serait intéressant d'en avoir une démonstration directe.

  \bigskip
  
  \noi 2.8. {\it Interprétation homologique de $h_1$ et $h_2$}.  
  . 
  \smallskip
 
  Soit $H_n(G,\Z), n \geqslant 0$, le $n$-ème groupe d'homologie de $G$ à coefficients
  dans $\Z$ (pour l'action triviale de $G$ sur $\Z$).   
   
\smallskip
\noi {\bf Théorème 2.9.} {\it Pour $n = 1$ et $n=2$, $H_n(G,\Z)$ est un $2$-groupe abélien élémentaire de rang $h_n$.}

     \smallskip
\noi {\it Démonstration}.

\smallskip
$\bullet \ n=1$. On a  $H_1(G,\Z) = G^{\ab} = G/D(G)$. \footnote{Si $G$ est un groupe, $D(G)$ désigne le groupe dérivé de $G$, autrement dit le sous-groupe engendré par les commutateurs; le groupe $G^{\ab}= G/D(G)$ est le plus grand quotient abélien de $G$.} Toute réflexion $s$ définit un élément de $G^{\ab}$ qui ne dépend que de la classe de conjugaison de $s$ ; d'où un homomorphisme $\theta_1: \F_2^{\Inv_1} \to G^{\ab}$. 

\noi (Rappelons que, si $X$ est un ensemble, $\F_2^X$ désigne le groupe des fonctions sur 
 $X$ à valeurs dans $\F_2$; lorsque $X$ est fini, ce groupe est souvent noté $\F_2[X]$.)
 
L'homomorphisme $\theta_1: \F_2^{\Inv_1} \to G^{\ab}$ 
est surjectif puisque $G$ est engendré par $\Inv_1$. En particulier, $G^{\ab}$ est un 2-groupe abélien élémentaire de rang $\leqslant h_1$. D'autre part, soit $\sigma \in \Inv_1$, et soit $S$ un ensemble de réflexions qui soit une base de $G$, au sens de la théorie des groupes de
Coxeter. Définissons une application $f: S \to  \F_2$ par $f(s) = 1$ si $s \in \sigma$ et $f(s) = 0$ sinon. Si $s,s' \in S$, et si $m(s,s')$ est l'ordre de $ss'$, on a
$m(s,s')(f(s)+f(s')) = 0$; en effet, c'est clair si $m(s,s')$ est pair; si $m(s,s')$ est impair,
$s$ et $s'$ sont conjugués, et l'on a $f(s)=f(s')$, d'où $f(s)+f(s')=0.$ Puisque  $G$
est défini par les relations  $(ss')^{m(s,s')}=~1$, cela montre que $f$ 
se prolonge en un homomorphisme de $G$ dans $\F_2$. Cet homomorphisme vaut $1$
sur les réflexions appartenant à $\sigma$, et $0$ sur les autres. Cela montre
que des réflexions deux à deux non conjuguées ont des images dans $G^{\ab}$ qui sont linéairement indépendantes. Cela entraîne que le rang de $G^{\ab}$ est $\geqslant h_1$, donc égal à $h_1$, et l'homomorphisme $\theta_1$ est bijectif.

\smallskip

$\bullet \ n=2$.
  En utilisant la formule de Künneth, on se 
  ramène au cas où $G$ est irréductible. D'après [IY 65], le groupe de cohomologie
  $H^2(G,\Q/\Z)$ est un 2-groupe abélien élémentaire dont le rang (noté $\kappa$ dans [IY 65])
  est donné explicitement pour chaque type $A_n,..., I_2(m)$. On constate (cas par cas) que $\kappa = h_2$. Comme $H^2(G,\Q/\Z)$ est dual de $H_2(G,\Z)$ (cf. [La 96], IV.5.5), on en déduit 2.9. 
  
\smallskip

\noi 2.10. {\it Compléments au th. $2.9$.}

\smallskip
    \noi  (a).  Dans le cas $n=1$ du th. 2.9, nous avons défini un isomorphisme  canonique $\theta_1:\F_2^{\Inv_1} \to H_1(G,\Z)$. On peut faire de même pour $n=2$,
    et expliciter un isomorphisme de $\F_2^{\Inv_2}$ sur $H_2(G,\Z).$
    Voici  comment. Si $u$ est une involution de degré 2, on a vu dans 2.3 que le
    groupe $D_u=G_u^{^-}$  est un groupe diédral  d'ordre divisible par 4. Un calcul standard montre que $H_2(D_u,\Z) \simeq \F_2$; d'où un  homomorphisme $f_u: \F_2 = H_2(D_u,\Z) \to H_2(G,\Z)$.
   Cet homomorphisme ne dépend que de l'image de $u$ dans $\Inv_2 $. La famille
   des $f_u$ définit donc un homomorphisme $\theta_2:\F_2^{\Inv_2} \to H_2(G,\Z)$.
   On peut démontrer que $\theta_2$ est un isomorphisme.

   \smallskip
   \noi (b). Au lieu d'utiliser l'homologie, on peut se servir de la cohomologie. C'est particulièrement commode  en dimension 2, puisque $H^2(G,A)$
   s'interprète comme le groupe des classes d'extensions centrales de $G$ par $A$. Le cas
   $A = \F_2$ peut se décrire explicitement de la manière suivante: 
   
     Soit $z \in H^2(G,\F_2)$. Si $u$ est une involution de $G$, définissons $f_z(u)\in \F_2$ comme 0 si la restriction de $z$ au sous-groupe $\langle u \rangle$ est $0$, et
     $f_z(u)=1$ sinon. Lorsque $u$ et $u'$ sont des involutions conjuguées, on a $f_z(u) = f_z(u')$. D'où une application, encore notée $f_z$, de $\Inv_G$ dans $\F_2$.
     On peut démontrer:
     
     \smallskip
     
     (b$_1$) {\it On a $z=0$ si et seulement si la restriction de $f_z$ à $\Inv_1 \cup \Inv_2$ est \ $0$.}
     
     \smallskip
     (b$_2$) {\it Toute application $f:\Inv_1 \cup \Inv_2 \to \F_2$ est de la forme $f_z$ pour un élément $z$ de $H^2(G,\F_2)$} (qui est unique d'après (b$_1$)).
   
     \smallskip
     
    \noi On a ainsi un {\it isomorphisme canonique} $H^2(G,\F_2) \simeq \F_2^{\Inv_1\cup \Inv_2}.$ 
     
 \noi  (Le fait que ces groupes soient isomorphes résulte du th. 2.9 par la formule des coefficients universels ([CE 56], VI.3.3a) - toutefois cette méthode ne donne pas un
     isomorphisme canonique.)

      \smallskip
   \noi (c). La condition $n=1,2$ du th. 2.9 est essentielle. Lorsque $n > 2$,
   il n'est pas vrai en général que $H_n(G,\Z)$ soit un 2-groupe,
ni qu'il puisse être engendré par $h_n$ éléments.

  \bigskip
{\bf $\S$3. Involutions maximales}
 
 \smallskip
\noi 3.1.   Soit $\rgr(G)$ le maximum des degrés des involutions de $G$; nous l'appellerons le {\it rang réduit} de $G$. Une
involution $u$ est dite {\it maximale} si son degré est égal à $\deg(u)=\rgr(G)$.

  Dans la terminologie de Springer ([Sp 74], $\S$ 4.2), cela signifie (si $G \neq 1$) que $u$ est un {\it élément régulier} d'ordre 2; Springer en donne plusieurs caractérisations,
  dont celle-ci: il existe un élément de $V^{^-}_u$ qui n'est fixé par aucun élément
  de $G$ distinct de 1; vu le th. 1.6, cela revient à dire qu'il n'existe aucune réflexion fixant  $V^{^-}_u$.
 
 \smallskip
 
 \noi {\bf Proposition 3.2}. {\it Soit $s$ une réflexion, et soit $G_s$ son centralisateur. Alors}:
 
  (a) {\it Toute involution maximale de $G$ est conjuguée à un élément de $G_s$.}
  
  (b) {\it Si $G$ est de type impair, toute involution de $G$ est conjuguée à un élément de $G_s$.}
  
   \smallskip

 \noi {\it Démonstration.}  On raisonne par récurrence sur l'ordre $|G|$ de $G$. On peut supposer que $G$ est irréductible. Si $-1$ appartient à $G$, c'est la seule involution maximale de $G$, et (a) est évident. On peut donc supposer
 que $-1 \not \in G$, donc que $G$ est de type impair, et il suffit de prouver (b).
 Or, si $u$ est une involution, il existe au moins une réflexion $s'$ qui commute à $u$; comme $s$ et $s'$ sont conjuguées (cf. 1.16) cela montre qu'un conjugué de $u$ commute à $s$, d'où (b).

 \smallskip

  Une conséquence de la prop.3.2.(a) est:
  
  \smallskip
  
  \noi {\bf Corollaire 3.3}.  \ $\rgr(G_s) = $ rgr$(G)$.

 \smallskip
\noi (Noter que $\rgr(G_s) = 1 + \rgr(G_s^{^+})$, avec les notations de 1.7.)
 
 \smallskip
 
   Autre conséquence:
   
   \smallskip
 
\noi {\bf Corollaire 3.4}. {\it Les involutions maximales sont conjuguées entre elles.}

  \smallskip
 
 \noi {\it Démonstration.}  On raisonne par récurrence sur  $|G|$. On peut supposer que $G$
 est irréductible de rang $>1$. Soit $s$ une réflexion. D'après la prop.3.2.(a), il suffit de prouver que les involutions maximales de $G_s$ sont conjuguées.
 Or $G_s$ n'est pas irréductible (cf. 1.7). L'hypothèse de récurrence s'applique donc
 à ses facteurs. D'où le résultat.
   
\noi (Pour une autre démonstration, de portée plus générale, voir [Sp 74], Theorem 4.2 (iv).)
 
  \medskip
  
 \noi 3.5. {\it Remarque}. Dans le cas où $G$ est irréductible et ne contient pas $-1$, le fait que toute involution maximale ait un conjugué dans
  $G_s$ permet une détermination de ces involutions par récurrence sur $|G|$. Cela s'applique aux cas suivants:
 
  $G$ de type $I_2(m), m $ impair \ : \ $G_s$ de type $A_1$; 
      
$G$ de   type $A_n, n \geqslant 2$ \ : \ $G_s$ de type $A_1 \times A_{n-2}$ (en convenant que $A_0=1$);
   
   $G$ de type $D_n, n $ impair \  : \ $G_s$ de type $A_1 \times D_{n-1}$;
   
 $G$ de  type $E_6$ \ : \ $G_s$ de type $A_1 \times A_5$.
   
   \smallskip
   
   \noi 3.6. {\it Agrandissement des involutions non maximales.}
   
   \smallskip
   
 \noi {\bf Théorème 3.7}. {\it    Soit $u$ une involution telle que  $\deg(u) = \rgr(G) - m$, avec $m >0$.
      Il existe $m$ réflexions $s_1,...,s_m$, commutant entre elles
     et commutant à $u$, telles que  $us_1\!\cdots\!s_m$ soit une involution maximale.}
     
     \smallskip
   \noi {\it Démonstration.}  Ici encore, on raisonne par récurrence sur $|G|$,
   que l'on peut supposer irréductible de rang $>1$. Soit $s$ une réflexion commutant à $u$.
   On a $u \in G_s$, et $G_s$ n'est pas irréductible. L'hypothèse de récurrence
   s'applique donc aux facteurs irréductibles de $G_s$, qui a le même rang réduit que $G$ d'après le corollaire 3.3.
   
   \smallskip
   
   \noi  3.8. {\it Remarque}. Le fait que $us_1\!\cdots\!s_m$ soit maximale entraîne que,
   pour tout $n \leqslant m$, l'involution $us_1\!\cdots\!s_n$ est de degré $\deg(u)+n$.

   \medskip

   \noi 3.9. {\it Involutions adjointes.}
   
     Soient  $u,v$  deux involutions. Disons que $v$ est {\it adjointe} de $u$ si
     les trois propriétés suivantes sont vérifiées:
     
     \smallskip 
     
    3.9.(a) \quad $u$ et $v$ commutent.
     
     3.9.(b) \quad $\deg(u)+\deg(v) = \deg(uv).$
     
     3.9.(c) \quad  $uv$ est une involution maximale, i.e. $\deg(uv) = \rgr(G)$.
     
     \smallskip
     
\noi     Noter qu'il y a symétrie: si $v$ est adjointe de $u$, alors $u$ est adjointe 
     de $v$.
     
     \medskip
     
     \noi {\bf Théorème 3.10}. {\it Toute involution a une adjointe, et cette adjointe est unique à conjugaison près.}
   
 \smallskip
 
 $\noi$ {\it Démonstration.}  Soit $u$ une involution et soit $m = \rgr(G)-\deg(u)$.
 Si $m=0$, $u$ est maximale, et sa seule adjointe est $1$. Si $m>0$, soient
 $s_1,...,s_m$ des réflexions ayant les propriétés du th. 3.7. L'involution $v= s_1\!\cdots\!s_m$ est adjointe de $u$. Cela démontre la partie ``existence'' du théorème.
 Pour l'unicité à conjugaison près, on raisonne par récurrence sur $|G|$. On remarque que, si $u$ et $v$ commutent, la propriété 
 3.9.(b) équivaut dire que $v$ fixe le sous-espace $V^{^-}_u$ de $V$. Si $v$ et $w$
 sont adjointes de $u$, elles appartiennent au fixateur $G^{^+}_u$ de $V^{^-}_u$, qui est un sous-groupe propre de $G$ (sauf si $u$=1, auquel cas le corollaire 3.4 montre que $v$ et $w$ sont conjuguées). On applique alors l'hypothèse de récurrence au $\cc$-sous-groupe $G^{^+}_u$.
 
 \smallskip
 
 \noi  3.11. {\it Application.} Soit $d = \rgr(G)$. Si $n$ est un entier $\geqslant 0$, rappelons (cf.  2.7) que $\Inv_n$ désigne
 l'ensemble des classes de conjugaison des involutions de degré~$n$. D'après le th. 3.10, la relation ``être adjointes'' définit une bijection entre $\Inv_n$ et $\Inv_{d-n}$; d'où $ |\Inv_n| =  |\Inv_{d-n}|$.
 
 \medskip
    
   \noi 3.12. {\it  Relations entre involutions maximales et chambres}.
   
   \smallskip
   
   Soit $\hh$ l'ensemble des hyperplans $H$ de $V$ tels qu'il existe une réflexion
   de $G$ fixant $H$. Soit $D = \bigcup_{H\in \hh} H$. Les composantes connexes
   de l'espace $V \sm D$ sont les {\it chambres} de $(V,G)$, cf. [Bo 68], $\S$V.1.3. Si $C$ et $C'$ sont deux chambres, il existe un unique $g\in G$
   tel que $gC = C'$ ([Bo 68],  $\S$V.3.2, th.1). En particulier, pour toute chambre $C$, il existe un unique $u_C \in G$ qui transforme $C$ en la chambre opposée $-C$; le carré de $u_C$ stabilise $C$, donc est égal à 1; autrement dit, $u_C$ est une involution.

\smallskip

\noi {\bf Théorème 3.13}. {\it Les $u_C$ sont les involutions maximales de $G$.}
   
   \smallskip
   \noi {\it Démonstration.} Si $-1 \in G$, on a $u_C=-1$ et c'est l'unique involution maximale de $G$. Si $-1 \not \in G$, on peut supposer $G$ irréductible. Le cas où $G$
   est un $I_2(m)$, $m$ impair, est immédiat; dans tout autre cas, $G$ est le groupe de Weyl d'un système de racines $R$ qui est irréductible et réduit. Pour prouver que $u_C$ est maximale, d'après le th. 3.7, il suffit de prouver qu'il n'existe aucune réflexion $s$ qui commute à $u_C$ et soit telle que $\deg(u_Cs) = \deg(u_C)+1$; cela équivaut à
   montrer que $u_C$ ne fixe aucun élément de $R$. Soit $\overline C°$ le cône
   dual de l'adhérence de $C$, autrement dit l'ensemble des $x\in V$ tels que $\langle x,y \rangle \geqslant 0$ pour tout $y\in C$. On sait que $R$ est contenu dans $\overline C° \cup -\overline C°$,
   cf. [Bo 68], $\S$VI.1.6. Comme $u_C$ échange $\overline C°$ et $-\overline C°$, et 
   que $\overline C° \cap -\overline C° = \{0\}$, cela entraîne que $u_C$ ne fixe aucun élément de~$R$. D'où le théorème.
   
   \medskip
   
   \noi 3.14. {\it Remarque.} Dans le cas où $G$ est irréductible, et ne contient pas $-1$, le th. 3.13, combiné aux Tables du Chap.VI de [Bo 68], donne une construction
   simple des involutions maximales: le graphe de Coxeter possède un seul automorphisme $\varepsilon$ d'ordre 2 (par exemple l'échange des deux ``bras'' pour le type $E_6$ ou le type $D_n, n $ impair); cet automorphisme définit un automorphisme
   externe $e$ de $(V,G)$; l'opposé $w_0=-e$ de $e$ appartient à $G$, et transforme
   la chambre canonique en son opposée, donc est une involution maximale d'après le th. 3.13. Cela montre que $\rgr(G)$ est égal à $\rg(G)- m/2$, où $m$ est le nombre
   des sommets du graphe non fixés par $\varepsilon$. Ainsi, pour $G$ de type $E_6$, on a
   $m=4,\  \rg(G)=6$, d'où $\rgr(G)=4$.
   
   \medskip
   
   \noi 3.15. {\it Centralisateur d'une involution maximale.}
   
   \medskip
   Soit $u$ une involution maximale de $G$.  Soit  $V = V_u^{^+} \oplus V_u^{^-}$  la décomposition de $V$ suivant les valeurs propres de $u$, et soit
   $G_u^{^-}$  (resp. $G_u^{^+}$) le fixateur de
   $V_u^{^+}$ (resp. $V_u^{^-}$), cf. 2.3. Puisque $u$ est maximale, le groupe $G_u^{^+}$
   ne contient aucune réflexion; comme c'est un $\cc$-sous-groupe de $G$, il est égal à $1$. Il en résulte que {\it le centralisateur $G_u$ de $u$  opère fidèlement} sur $V_u^{^-}$.

   \medskip
   
  \noi {\bf Théorème 3.16}. ([Sp 74], Theorem 4.2 (iv)). {\it  Le couple $(V_u^{^-},G_u)$ est un couple de Coxeter.}
  
  \noi En particulier, $G_u$ est un groupe de Coxeter, bien que n'étant pas en général un $\cc$-sous-groupe de $G$. En fait, comme la démonstration le montrera,
  $G_u$ est engendré par des involutions de degré 1 ou 2 de $G$.
  
  \noi {\it Remarque}. On peut se demander si cet énoncé reste valable pour des involutions quelconques; G. Röhrle m'a signalé que c'est le cas; je reviendrai sur cette question dans ``{\it Groupes de Coxeters finis} : {\it centralisateurs d'involutions} '', à paraître.
  
  \smallskip
  
  \noi {\it Démonstration}. Il suffit de traiter le cas où $G$ est irréductible et où
  $u \neq -1$, ce qui signifie que $G$ est de l'un des types $I_2(m), m$ impair,
  $A_n, n \geqslant 3$, $D_n, n$ impair, et $E_6$. On a alors:
    
 \medskip 
    type $I_2(m), m $ impair  $\Longrightarrow$  $G_u$ est de type $A_1$.
    
    \smallskip
   
    type $A_n, n \geqslant  3$ \hspace{9mm} $\Longrightarrow$   $G_u$ est de type $B_{[(n+1)/2]}.$
    
    \smallskip
    
    type $D_n, n$ impair \hspace{4mm}  $\Longrightarrow$  $G_u$ est de type $B_{n-1}$.
    
    \smallskip
    
   type $E_6$ \hspace{20mm}  $\Longrightarrow$  $G_u$ est de type $F_4$.
    
    \medskip
 \noi   De plus $G_u$ est engendré par les réflexions qu'il contient (celles qui commutent  à $u$), et les involutions de degré 2 de la forme  $\sigma = ss'$, où
    $s$ et $s'$ sont des réflexions distinctes, commutant entre elles, et échangées
    par conjugaison par $u$ (i.e. $s' = usu$); une telle involution est de degré 2, vue comme élément de $G$, mais de degré 1, vue comme élément de $G_u$ opérant sur $V_u^{^-}$.
    
    \smallskip
    
     \noi Ces assertions se vérifient par calcul explicite:

      \smallskip
      Si $G$ est de type $I_2(m)$,$m$ impair, on a $G_u = \langle u\rangle$. 
      
      Si $G$ est de type $A_n$, donc isomorphe à $\Sym_{n+1}$,
      on peut prendre pour $u$ le produit des $m$ transpositions (12), (34),...,
      où $m = [\frac{n+1}{2}]$; le groupe $G_u$ est engendré par ces transpositions, et par le groupe $\Sym_m$ qui les permute; c'est un groupe de type $B_m$.
      
      Si $G$ est de type $D_n, n$ impair, on prend pour $u$ l'involution qui fixe la dernière variable $x_n$, et change le signe de $x_1,...,x_{n-1}$; le groupe $G_u$
      est le groupe $B_{n-1}$, opérant de façon naturelle sur $x_1,...,x_{n-1}$, et
      opérant sur $x_n$ par un signe convenable (de façon à obtenir un plongement
      de $B_{n-1}$ dans $D_n$).

     Si $G$ est de type $E_6$, on prend $u = - e$ , où  $e$  est l'automorphisme externe qui permute $\alpha_1$ et $\alpha_6$, ainsi que $\alpha_3$ et $\alpha_5$, cf. figure ci-dessous:

  \vspace{6mm}
  
\hspace{25mm}$\alpha_{_1} \hspace{6mm}\alpha_{_3}\hspace{5mm}\alpha_{_4}\hspace{6mm} \alpha_{_5}\hspace{5mm} \alpha_{_6}$

\hspace{25mm}$\circ$-------$\circ$------$\circ$-------$\circ$------$\circ$

\hspace{43mm} $\mid$

\hspace{43mm} $\circ \ \alpha_{_2}$

\medskip
    
      Soient $\beta_1=\alpha_1+\alpha_6$ et $\beta_3=\alpha_3+\alpha_5$.
      Soit $S = \{\beta_1, \beta_3, \alpha_4,\alpha_2\}$; c'est une base de
      $V_u^{^-}$.
  
      Si le produit scalaire sur $V$ est normalisé de telle sorte que $\alpha_i.\alpha_i=1$ pour  $i=1,...,4$, on a $\beta_i.\beta_i=2$ pour  $i=1,2$, et:
       $$ \beta_1.\beta_3 = -1, \ \ \beta_3.\alpha_4=-1, \ \ \alpha_4.\alpha_2=- 1/2;$$
       les autres produits scalaires sont nuls. Il en résulte que  $S$ est la base d'un système de racines $R$ de type $F_4$ dans  $V_u^{^-}$, les racines longues étant
       $\beta_1$ et $\beta_2$:
      
      \medskip
     
\hspace{25mm}$\beta_{_1} \hspace{6mm}\beta_{_3}\hspace{6mm}\alpha_{_4}\hspace{7mm}\alpha{_2} $

\hspace{25mm}$\bullet$-------$\bullet$=$>$=$\circ$-------$\circ$

\medskip

On obtient ainsi un sous-groupe de Coxeter $G'$  de $G_u$ de type $F_4$. Ce sous-groupe est égal à $G_u$. En effet tout élément de $G_u$ définit un automorphisme
du système de racines $R$ (car $R$ est l'ensemble des racines de $G$ invariantes par $u$), donc appartient à  $G'$, cf. [Bo 68], $\S$VI.4.9.

\medskip

\noi 3.17. {\it Relations entre involutions maximales et invariants polynomiaux.}
   
 \smallskip
 
 Supposons que $(V,G)$ soit réduit; soit $n = \rg G$. Soit $S$ l'algèbre des fonctions polynômes sur $V$ qui sont $G$-invariantes;  c'est une algèbre graduée de polynômes [Bo 68], $\S$V.5.3, th.3), dont le corps des fractions est de degré de transcendance $n$ sur $\R$. Soient  $d_1,...,d_n$ les degrés caractéristiques de $S$ ; on a $\prod d_i = |G|$, cf. [Bo 68], $\S$V.5.3, cor. au th.3.
   
 \medskip
 
 \noi {\bf Théorème 3.18.} {\it Soit $u$ une involution maximale de $G$, et soit $G_u$ son centralisateur}.

   (a) {\it Le nombre des $i$ tels que $d_i$ soit pair est égal à $\rgr(G)$.} 
  
  (b) {\it Les $d_i$ pairs sont les mêmes que ceux du groupe de Coxeter $G_u$.}
  
  (c) {\it Le produit des $d_i$ impairs est égal au nombre des involutions maximales
  de $G$.}
  
  \medskip
  
 \noi {\it Démonstration}. Les assertions (a) et (b) sont des cas particuliers du
 th. 4.2 (iii) de [Sp 74]. On déduit de (b) que le produit des $d_i$ pairs
 est égal à $|G_u|$; comme le produit des $d_i$ est égal à $|G|$,
 cela montre que produit des $d_i$ impairs est égal à $(G:G_u)$; d'où (c).
  
   \medskip
  
  \noi {\bf Corollaire 3.19.} {\it Le nombre des involutions maximales de $G$ est impair.}
  
\vspace{11mm}
 
 {\bf $\S$4. Cubes}. 
 
 \medskip

     \noi 4.1. Un {\it cube} de $G$ est un  $\cc$-sous-groupe abélien $C$ de $G$ (cf. [Se 18]), autrement dit un sous-groupe engendré par des réflexions commutant deux à deux;  c'est un groupe de Coxeter de type $A_1 \times \cdots \times A_1$. Le nombre de facteurs $A_1$ est le rang de $C$.
      
      \smallskip
     
     Soit $d = \rg(C)$. On a $|C|  = 2^d$. L'ensemble $S$ des réflexions de $C$ 
     a $d$ éléments; on l'appelle la {\it base} de $C$. Tout élément de $C$ s'écrit de façon unique comme produit  $s_I = \prod_{s\in I} s$, où $I$ est une partie de $S$; on a
 $ \deg(s_I) = |I|.$ Le cube $C$ contient un et un seul élément de degré $d$: celui qui correspond à $I=S$; nous appellerons cet élément {\it l'extrémité} de $C$. {\it Toute involution est extrémité d'au moins un cube}, d'après la prop.2.5.
 
 \smallskip
 
 \noi 4.2. {\it Centralisateur d'un cube.}
 
 \smallskip
   Soit $C$ un cube de rang $d$, soit $S$ sa base et soit $u$ son extrémité.
   Soit $V=  V_u^{^+} \oplus V_u^{^-}$  la décomposition de $V$ associée à $u$, cf.
   2.2. L'espace $V_u^{^+}$ est l'ensemble des points de $V$ fixés par $C$. L'espace $V_u^{^-}$ est somme directe des droites $D_s = V_s^{^-}$, pour $s \in S$. Soit $G_u^{^+}$ le fixateur de $V_u^{^-}$, et soit $G_C$ le centralisateur de $C$ dans $G$. On a:
   
   \smallskip
   
   \noi {\bf Proposition 4.3}.  $G_C = C \times G_u^{^+}$.
   
   \smallskip
   
   \noi {\it Démonstration.} L'inclusion $C \times G_u^{^+} \subset G_C$  est évidente. Inversement, tout $g \in G_C$ stabilise $V_u^{^-}$ ainsi que les droites
   $D_s$; sur chacune de ces droites il agit, soit trivialement, soit par changement de signe. En le multipliant par un élément convenable  $c$ de $C$, on obtient un
   élément $gc$ qui fixe chaque $D_s$, donc  $V_u^{^-} $; d'où $gc \in G_u^{^+}$.

   \smallskip
   
   \noi {\bf Corollaire 4.4}. {\it Le centralisateur $G_C$ de $C$ est un $\cc$-sous-groupe de $G$.}

   \medskip
   
   \noi {\bf Théorème 4.5}. {\it Supposons $G$ de type impair.  Soit $u$ une involution de $G$. Si $\deg(u) \leqslant 2$, il n'existe qu'un seul cube d'extrémité $u$. Il en est de même si $\deg(u)=3$ et si $G$ est cristallographique.
   Si $\deg(u) =4$, et si $G$ est cristallographique, le nombre des cubes d'extrémité $u$ est, soit $1$, soit $3$.}
   
   \smallskip
   
   \noi {\it Démonstration.}  Soit $G_u^{^-}$ le fixateur de $V_u^{^+}.$ Les cubes d'extrémité $u$ sont les cubes maximaux de $G_u^{^-}$. Le couple $(V_u^{^-}, G_u^{^-})$ est un couple de Coxeter
de type impair, et $G_u^{^-}$ contient $-1$; son rang est $d=\deg(u)= \dim V_u^{^-}$. Lorsque $d \leqslant 2$, cela entraîne que $G_u^{^-}$ est un produit de $d$ groupes de type $A_1$, autrement dit est un cube, et c'est le seul cube 
d'extrémité $u$. Si $G$ est cristallographique, il en est de même pour $d=3$, et si $d=4$, le groupe $G_u^{^-}$ est, soit de type $A_1\times A_1 \times A_1\times A_1$, soit de type  $D_4$; dans ce dernier cas,
le nombre de cubes d'extrémité $u$ est $3$.

\smallskip

\noi {\it Remarque.} L'hypothèse que $G$ est de type impair est essentielle:
 si $G$ est de type $B_2$ (resp. $G_2$), il y a 2 (resp. 3) cubes
d'extrémité $-1$. Celle que  $G$  est cristallographique est aussi essentielle: 
si $G$ est de type $H_3$ (resp. $H_4$), il y a 5 (resp. 75) cubes
d'extrémité $-1$.

   \smallskip
   
   \noi 4.6. {\it Cubes maximaux.}
   
   \smallskip
    Un cube est dit {\it maximal} si son extrémité est une involution maximale,
    autrement dit si son rang est égal à $\rgr(G)$. Cette terminologie est justifiée par le résultat suivant:
    
    \smallskip
    \noi {\bf Proposition 4.7}. {\it Tout cube est contenu dans un cube maximal.}
   
   \smallskip
   
   \noi {\it Démonstration}. Soit $C$ un cube; soit $u$ son d'extrémité et soit $S$ sa base. Soit $m= \rgr(G)-\deg(u)$. Si $m=0$, $C$ est maximal. Supposons
   $m > 0$. D'après le th.3.7, il existe $m$ réflexions
   $s_1,...,s_m$, commutant entre elles et commutant à $u$, telles que $us_1...s_m$ soit une involution maximale. On a $\deg(us_i) = \deg(u)+1$ pour tout $i$. Cela entraîne que les droites $V_{s_i}^{^-}$ sont orthogonales à $V_u^{^-}$, donc que les $s_i$ commutent aux éléments de $S$. Le cube de base $S \cup \{s_1,...,s_m\}$ est un cube maximal
   contenant~$C$.   
   
   \medskip
   
   \noi {\bf Proposition 4.8}. {\it Un cube maximal est égal à son centralisateur.}
   
    \smallskip
   
   \noi {\it Démonstration}. Soit $C$ un cube maximal et soit $u$ son extrémité.
   D'après la prop.4.3, le centralisateur de $C$ est $C \times G_u^{^+}$, et d'après
   3.15, on a $G_u^{^+}=1$.
   
   \medskip
   
   \noi 4.9. {\it Les cubes d'un groupe de type impair.}
   
   \smallskip
   Les cubes d'un groupe de Coxeter de type impair ont des propriétés analogues
   à celles des $2$-sous-groupes de Sylow d'un groupe fini. Tout d'abord:
   
   \smallskip
   
   \noi {\bf Théorème 4.10}. {\it Supposons $G$ de type impair. Soit $u$ une involution. Soient $C$ et $C'$ deux cubes d'extrémité $u$. Alors $C$ et $C'$ sont conjugués.}
   
   (Noter que, si $g\in G$ est tel que $gCg^{-1} = C'$, on a $gug^{-1} = u,$ donc  $g$
   appartient au centralisateur $G_u$ de $u$.)
   
   \smallskip
  
   \noi {\it Démonstration.} On raisonne par récurrence sur $|G|$, et l'on peut supposer $G$ irréductible et d'ordre $>2$. Les cubes $C$ et $C'$ sont contenus dans
   $G_u^{^-}$. Si $G_u^{^-} \neq G$, l'hypothèse de récurrence montre que $C$ et $C'$ sont conjugués dans $G_u^{^-}$, donc dans $G$. Si $G_u^{^-} = G$, on a
  $V_u^{^+} = 0$, d'où $u=-1$. Soient $s, s'$ deux réflexions appartenant respectivement à $C$ et à $C'$. Puisque $G$ est de type impair, il existe $g\in G$ tel que $s'=gsg^{-1}$, cf. th.1.15. Les cubes $C'$ et $gCg^{-1 }$ contiennent tous deux $s'$, donc sont contenus dans le centralisateur $G_{s'}$ de $s'$, et ont la même extrémité, à savoir $-1$. Comme les facteurs irréductibles
  de $G_{s'}$ sont d'ordre $ < |G|$, l'hypothèse de récurrence montre
  que $C'$ et $gCg^{-1 }$ sont conjugués dans $G_{s'}$, d'où le fait que $C$ et $C'$
  sont conjugués dans $G$.
  
  \smallskip
  
   \noi {\bf Corollaire 4.11}. {\it Si $G$ est de type impair, les cubes maximaux de $G$ sont conjugués entre eux.}
   
       \smallskip
Cela résulte de 3.4 et 4.10.

\smallskip

\noi {\it Remarque.} Supposons que $G$ soit irréductible, mais pas de type impair.
Alors les cubes maximaux de $G$ sont conjugués si $G$ est de type $I_2(2m), m $ impair (par exemple de type $G_2$), mais ils ne le sont pas dans les autres cas:
$I_2(2m), m$ pair, $B_n$ et $ F_4$.
  
  \smallskip
  
  \noi {\bf Proposition 4.12}. {\it Si $G$ est de type impair, le nombre des cubes d'extrémité une involution donnée est impair.} 
  
  \smallskip
   
   \noi {\it Démonstration.} Le même argument que ci-dessus permet de se ramener
   au cas où $G$ est irréductible et où l'involution est égale à $-1$. On vérifie cas par cas (voir $\S$5) que ce nombre est impair: 
   
   type \ $A_1$\  $\Longrightarrow$ 1;
   
     ------ \ $D_n, n$ pair $\Longrightarrow$  $1\!\cdot\!3\!\cdot\!5\cdots\!(n-1)$;
     
       ------ \ $E_7$ $\Longrightarrow$ $3^3\!\cdot\! 5$;
     
       ------ \ $E_8$ $\Longrightarrow$ $3^4\!\cdot\! 5^2$;
     
       ------ \ $H_3$ $\Longrightarrow$ $5$;
     
       ------ \ $H_4$  $\Longrightarrow$ $3\!\cdot\! 5^2$;
     
       ------ \ $I_2(m), m$ impair $\Longrightarrow$ $m$.
       
       \smallskip
             
       \noi {\bf Corollaire 4.13}. {\it Si $G$ est de type impair, le nombre des cubes maximaux
       de $G$ est impair.}
       
       \smallskip
           En effet, ce nombre est égal au produit du
        nombre des involutions maximales par le nombre des cubes maximaux d'extrémité donnée; or les deux facteurs sont impairs d'après 3.18 et 4.12.
     
     \smallskip
      \noi {\bf Corollaire 4.15.} {\it Si $G$ est de type impair, tout $2$-sous-groupe de Sylow de $G$ contient un sous-groupe normal qui est un cube maximal}. 
        
        \smallskip
     {\it Démonstration.}    Soit $P$ un  $2$-sous-groupe de Sylow de $G$. Faisons-le opérer par conjugaison sur l'ensemble $X$ des cubes maximaux. Comme $|X|$  est impair, l'une des orbites de $P$ a un seul élément. Il existe donc un cube maximal $C$
     qui est normalisé par $P$. Le produit $P\!\cdot\! C$ est un $2$-groupe contenant $P$, donc est égal à $P$. Cela montre que $C$ est un sous-groupe normal de $P$.

\medskip

\noi 4.15. {\it Contrôle de la fusion dans un cube maximal}.

\smallskip  

L'énoncé ci-dessous est analogue au classique théorème de Burnside
sur la fusion dans le centre d'un groupe de Sylow:

\smallskip

\noi {\bf Théorème 4.16}. {\it Supposons $G$ de type impair. Soit $C$ un cube maximal de $G$ et soit $N_C$ son normalisateur. Soient $A$ et $B$ deux parties de $C$ et soit $g\in G$ tel que $gAg^{-1} = B$. Il existe alors $n \in N_C$ tel que
$gag^{-1}=nan^{-1}$ pour tout $a\in A$.}

\noi [En particulier, $A$ et $B$ sont conjuguées dans $N_C$.]
   
   \smallskip
   
   \noi {\it Démonstration} (imitant celle du théorème de Burnside).
    Soit $G_A$ l'intersection des centralisateurs des éléments de $A$; définissons
    $G_B$ de manière analogue. On a $C \subset G_A \cap G_B$ et $gG_Ag^{-1}=G_B$. Les deux cubes $C$ et $gCg^{-1}$ sont contenus dans $G_B$; soit $H =
    \langle C,gCg^{-1}\rangle$ le sous-groupe de $G_B$ qu'ils engendrent; c'est un $\cc$-sous-groupe de $G$. D'après le corollaire 4.11, appliqué à $H$, il existe $h\in H$ tel que
    $hCh^{-1} = gCg^{-1}$, d'où $h^{-1}g \in N_C$. Soit $n = h^{-1}g$.
     Si $a\in A$, on a $nan^{-1} =h^{-1}\!\cdot\! gag^{-1}\!\cdot\! h = gag^{-1}$
    puisque $gag^{-1}$ appartient à $B$ et $h$ centralise $B$. 
     
   \medskip
   
   \noi 4.17. {\it Description des classes d'involution d'un groupe de type impair à partir d'un cube maximal.}
   
   \smallskip
    Conservons les notations et hypothèses du th.4.16. Posons $\Phi_C= N_C/C$.
    Le groupe $\Phi_C$ agit par conjugaison sur $C$, d'où un homomorphisme $\Phi_C \to \Aut(C)$, qui est injectif puisque $C$ est son propre
    centralisateur (4.8). On peut donc identifier $\Phi_C$ à un sous-groupe
    de $\Aut(C)$. Soit $\Inv_G$ l'ensemble des classes de conjugaison d'involutions
    de $G$. L'application naturelle  $C \to \Inv_G$ est constante sur chaque orbite
    de $\Phi_C$; elle définit donc une application de l'ensemble quotient $C/\Phi_C$ dans $\Inv_G$.
    
    \smallskip
    \noi {\bf Proposition 4.18}.  {\it L'application $C/\Phi_C \to \Inv_G$ est bijective.}
    
    \smallskip
    \noi {\it Démonstration.} Si $u$ est une involution, elle est contenue dans un cube; ce cube est contenu dans un cube maximal, lequel est conjugué de $C$ d'après le cor.4.11; cela montre que $C/\Phi_C \to \Inv_G$ est surjectif. L'injectivité résulte du th.4.10.
    
    \medskip
    On peut reformuler l'énoncé ci-dessus de façon plus explicite.  Soit $S$ la base
    de $C$ et soit $n$ un entier $\geqslant 0$. Le  groupe $\Phi_C$ s'identifie à {\it un groupe de permutations de $S$}. Il opère donc sur l'ensemble $S_n$ des parties à $n$ éléments de $S$; l'application $I \mapsto s_I$ de 4.1 donne une bijection de $S_n$ sur l'ensemble des éléments de $C$ de degré $n$; d'où
    une application $S_n \to \Inv_n$  (rappelons que $\Inv_n$ désigne l'ensemble des classes de conjugaison des involutions de degré $n$, cf. 2.7); d'après la prop.4.18,
    deux éléments $I, I'$ de $S_n$ ont même image dans $\Inv_n$ si et seulement
    si ils sont dans la même  $\Phi_C$-orbite. En résumé:
    
       \smallskip 
    \noi {\bf Proposition 4.19}. {\it Les éléments de $\Inv_n$ correspondent bijectivement aux orbites
    de $\Phi_C$ opérant sur l'ensemble des parties à $n$ éléments de la base de $C$.}
    
    \noi (Rappelons que $G$ est supposé de type impair.)
    
       \smallskip 
       
       Ainsi, les classes de conjugaison des
       involutions de $G$ se lisent sur l'action de $\Phi_C$ sur $S$; par exemple, l'opération d'adjonction de 3.9 transforme une partie de $S$ en son complémentaire. On en verra d'autres exemples dans 5.5 et 5.6  lorsque $G$ est de type $E_7$ ou $E_8$.

   \bigskip

 {\bf $\S$5. Description des différents types irréductibles}. 
 
 \smallskip
 
 Cette section rassemble quelques résultats sur les divers types de groupes de Coxeter finis irréductibles, leurs involutions et leurs cubes.
 
 Les notations sont celles des $\S\S$ 1,...,4, et aussi celles du chap. VI de [Bo 68];
 pour alléger l'écriture, on écrit ``$X$''  pour ``un groupe de Coxeter de type $X$''.
 
  Si $s$ est une réflexion, on note
 $G_s^{^+}$ le fixateur de la droite $V^{^-}_s$ correspondante; le centralisateur de $s$
 est $\{1,s\} \times G_s^{^+}$, cf. prop.1.8.

 \medskip
 
 \noi 5.1. {\it Type $A_{n-1}$ - ordre $n!$}
 
 \smallskip
 
   Soit $n > 1$. On a $A_{n-1} =\Sym_n$,
   agissant sur  $V = \R^n$ par permutation des coordonnées. Les réflexions sont les transpositions; les groupes $G_s^{^+}$ correspondant sont de type $A_{n-3}$ (en convenant
   que $A_m = 1$ si $m \leqslant 0$). 
   
   On a  $\rgr(A_{n-1})= [n/2].$ 
   
   Pour tout  $d \leqslant [n/2]$ il y a une seule classe d'involutions de degré $d$: les produits de $d$  transpositions à supports disjoints; chacune de ces involutions est l'extrémité d'un seul cube. Le cas $d = [n/2]$ est celui des  involutions maximales; leur nombre est égal au produit $f(n) = 1\!\cdot\!3\!\cdot\! 5\!\cdots$ des entiers impairs $\leqslant n$.
   
    Le centralisateur d'une involution maximale est un groupe de Coxeter de type $B_{[n/2]}$, cf. 3.16; le groupe $\Phi$ associé à un cube maximal (au sens de 4.17) est isomorphe à $\Sym_{[n/2]}$.
      
   \medskip
 
 \noi 5.2. {\it Type $B_n$ - ordre $2^nn!$}
  
 \smallskip 
 
   Soit $n\geqslant 2$. Le groupe  $B_n$  est le groupe des permutations et changements de signe de $n$ variables; c'est un sous-groupe
   du groupe orthogonal $\O_n(\R)$. Soit $(e_1,...,e_n)$ la base canonique de $\R^n$; l'action de $B_n$ sur les $n$ ensembles $\{e_1,-e_1\}$, ..., $\{e_n,-e_n\}$
définit une surjection $B_n \to \Sym_n$,
   dont le noyau est un cube de rang $n$, que nous noterons ${\bf2}^n$ (sous-groupe ``diagonal''). On a donc $B_n={\bf2}^n\!\cdot\!\Sym_n$ (produit semi-direct).
   
  Le groupe $B_2$ est isomorphe au groupe diédral $I_2(4)$ d'ordre 8.

   Il y a deux types de réflexions: celles, dites {\it longues}, dont l'image dans $\Sym_n$ est une transposition, et celles, dites {\it courtes} qui appartiennent à $I_n$. Les premières correspondent aux racines $\alpha= \pm e_i \pm e_j$ telles que $\alpha\!\cdot\!\alpha=2$; elles permutent un couple $(e_i,e_j)$ ou $(e_i,-e_j)$, $i \neq j$. Leur nombre est $n(n-1)$. Les groupes  $G_s^{^+}$ correspondants sont de type $A_1\times B_{n-2}$ (en convenant que $B_1 = A_1$ et $B_0= 1$). Les réflexions du second type
   correspondent aux racines $ \alpha= \pm e_i$ telles que $\alpha\!\cdot\!\alpha=1$; elles changent de signe un et un seul
 des $e_i$. Leur nombre est $n$. Elles forment une base du cube ${\bf2}^n$ défini
   ci-dessus. Les groupes  $G_s^{^+}$ correspondants sont de type $B_{n-1}$ (en convenant que $B_1=A_1$).
      
   Une involution $u$ de $B_n$ a deux invariants $a$ et $b$: le nombre $a$ des $e_i$ qui sont transformés en leur opposé, et le degré $b$ de l'image de $u$ dans
$\Sym_n$. La somme $a+b$ est égale au degré $d$ de $u$. A
  conjugaison près, et à des changements de signe près des $e_i$,
 $u$ transforme $e_1,...,e_a$ en $-e_1,...,-e_a$ et permute les $b$ couples
 $(e_{a+1},e_{a+2}), (e_{a+3},e_{a+4}), ...,(e_{a+2b-1},e_{a+2b})$. Les classes de conjugaison d'involutions de degré  $d$  correspondent ainsi aux couples ($a,b)$ d'entiers $\geqslant 0$ tels que $a+b=d$ et $a+2b \leqslant n$. Une involution adjointe (au sens de 3.9) à une involution d'invariants $(a,b)$ a pour invariants $(n-a-2b,b)$.
 L'unique involution maximale est $-1$; ses invariants sont $(n,0)$.
 
 Soit $u$ une involution d'invariants $(a,b)$. Les cubes d'extrémité $u$ ont un invariant supplémentaire: un entier  $c$
 tel que $0 \leqslant c \leqslant a/2$, caractérisé par le fait que la base du cube contient $a-2c$ réflexions courtes et $b+2c$ réflexions longues (autre caractérisation:
 l'image du cube dans $\Sym_n$ est un cube de rang $b+c$). Le cas
 $c=0$ correspond à un unique cube: celui dont la base est évidente d'après la description
 ci-dessus. Le cas $c =1$ se déduit du précédent en remplaçant deux réflexions
 courtes, par exemple $e_1 \leftrightarrow -e_1$ et $e_2 \leftrightarrow -e_2$  par les deux réflexions longues $e_1 \leftrightarrow e_2$ et $e_1 \leftrightarrow -e_2$. On procède de même pour $c = 2, ...,[a/2]$. Pour $c$ fixé, les cubes correspondants sont conjugués et leur nombre est égal à celui des involutions de $\Sym_a$ de degré $c$, autrement dit  $a\choose 2c$$f(a)$, où $f(a)$ est défini comme dans 5.1; c'est un nombre impair si $a=2c$ ou $a=2c+1$. Un cube maximal a pour invariants $(a,b,c) = (n,0,c)$ avec $c\leqslant n/2$; pour $c=0$, il y a un seul
 tel cube, à savoir ${\bf2}^n$.
 
  Lorsque $n$ est impair, $B_n =\{\pm1\} \times D_n$, où $D_n$ est le sous-groupe
  de $B_n$ engendré par les réflexions longues, cf. 5.3. Par exemple, $B_3$ est
  isomorphe à $\{\pm1\} \times D_3 \simeq \{\pm1\} \times A_3 \simeq \{\pm1\} \times \Sym_4$.
 
  \medskip
 
 \noi 5.3. {\it Type $D_n$ - ordre $2^{n-1}n!$}
 
 \smallskip
   Soit $n \geqslant 3$. Le groupe $D_n$ est le sous-groupe d'indice 2 de $B_n$ formé des permutations et changements de signes en nombre pair. Il est souvent noté ${\bf2}^{n-1}\!\cdot\!\Sym_n$. On a $\rgr(D_n) = n$ (resp. $n-1$) si $n$ est pair (resp. impair).
   
Les réflexions de $D_n$ sont les réflexions longues de $B_n$; les groupes $G_s^{^+}$ correspondants sont de type $A_1 \times D_{n-2}$ (en convenant que $D_1 = 1$ et $D_2 = A_1 \times A_1$). Les involutions sont celles de $B_n$ dont l'invariant $a$ est pair. Deux involutions de $D_n$ sont conjuguées dans $D_n$ si et seulement si elles le sont dans $B_n$ (autrement dit si elles ont les mêmes invariants $a$
   et $b$), {\it à une exception près}: quand $n$ est pair, la classe de conjugaison
   de $B_n$ d'invariants $(0,n/2)$ se décompose dans $D_n$ en deux classes distinctes. Par exemple, pour $n=6$ ces classes sont représentées par:
   
     $(e_1 \leftrightarrow e_2)\cdot(e_3\leftrightarrow e_4)\cdot(e_5 \leftrightarrow e_6)$ \ \ et \ \ $(e_1 \leftrightarrow e_2)\cdot(e_3\leftrightarrow e_4)\cdot(e_5 \leftrightarrow -e_6)$; 
     
  \noi   c'est la parité du nombre de signes ``$-$'' qui caractérise
     la classe de conjugaison.
     
     Pour $n=4$, ces deux classes, et celle d'invariants $(2,0)$, sont permutées par les automorphismes extérieurs de G (trialité). Avec les notations de [Bo 68], $\S$VI.4.8, on peut les représenter par $s_{\alpha_1}s_{\alpha_3}$, $s_{\alpha_1}s_{\alpha_4}$ et $s_{\alpha_3}s_{\alpha_4}$.

    L'image par la surjection $D_n \to \Sym_n$ d'un cube maximal de $D_n$ est un cube maximal de $\Sym_n \simeq A_{n-1}$. On obtient ainsi une bijection entre les cubes maximaux de $D_n$ et ceux de $A_{n-1}$.
     
       Pour $n$ pair, $-1$ appartient à $D_n$ et les cubes maximaux  sont ceux
     de $B_n$ d'invariants $(a,b,c) = (n,0,n/2)$; les groupes  $\Phi$ correspondants sont des groupes de Coxeter de type $D_{n/2}$. Pour
 $n$ impair, il y a $n$ involutions maximales; ce sont les involutions de $B_n$ d'invariants $(n-1,0)$; leurs centralisateurs sont de type $B_{n-1}$, cf. 3.16. Les cubes maximaux sont les cubes de $B_n$ d'invariants $(n-1,0,(n-1)/2)$ et les groupes $\Phi$ correspondants sont des groupes de Coxeter de type $B_{(n-1)/2}$. 
     
     \medskip
     \noi 5.4. {\it Type $E_6$ - ordre $2^7\!\cdot\!3^4\!\cdot\!5$ - } [ATLAS], p.26.
     
          \smallskip
          Le nombre des réflexions est $(78-6)/2=36$;  les groupes $G_s^{^+}$ correspondants sont de type $ A_5$, donc isomorphes à $\Sym_6$. Il en résulte
          que le rang relatif de $E_6$ est 4, et qu'il y a une seule classe d'involutions pour chaque degré
$\leqslant 4$.  Il y a 45 involutions maximales; leurs centralisateurs sont des groupes de Coxeter de type $F_4$, cf. 3.16. Chacune est l'extrémité de 3 cubes maximaux;  le groupe $\Phi$ associé à un tel cube est $\Sym_4$.

  Soit $R$ le réseau des racines et soit $P$ le réseau des poids. Le $\F_3$-espace vectoriel $V_5=R/3P$ est de dimension 5, et le produit scalaire de $V$ définit sur $V_5$ une forme bilinéaire symétrique $a$ non dégénérée. Si $g\in E_6$, l'action de $g$ sur $V_5$ définit un élément  $g_5$ du groupe orthogonal 
$\O(V_5,a)$; l'application $g \mapsto \det(g)g_5$ donne un isomorphisme 
 $E_6 \  \to \  \SO(V_5,a) \simeq \  \SO_5(\F_3).$

      Le groupe $E_6$ contient un $\cc$-sous-groupe de type $D_5$ qui est d'indice 27, à savoir le fixateur du poids fondamental $\omega_1$; l'espace homogène $E_6/D_5$ peut s'identifier à l'ensemble des 27 droites d'une surface
      cubique, cf. par exemple [Ma 74], $\S$26 et [BS 20].
      
      \medskip

     \noi 5.5. {\it Type $E_7$ - ordre $2^{10}\!\cdot\!3^4\!\cdot\!5\!\cdot\!7$} - [ATLAS], p. 46.
     
     \smallskip
     
     Le nombre des réflexions est $(133-7)/2=63$;
  les groupes $G_s^{^+}$ correspondants sont de type $D_6$. 
     
     Soit $R$ le réseau des racines et soit $P$ le réseau des poids. Le $\F_2$-espace vectoriel $V_6=R/2P$ est de dimension 6, et le produit scalaire de $V$ définit sur $V_6$ une forme bilinéaire alternée $a$ non dégénérée. L'action de $E_7$ sur $V_6$
     donne un isomorphisme: \  $E_7/\{\pm1\}  \ \to  \ \Sp(V_6,a) \ \simeq \  \Sp_6(\F_2)$. 
     
     Toute réflexion $s=s_\alpha$ définit un élément $\tilde{s}$ de $V_6\sm\{0\}$, à savoir l'image de $\alpha$ dans $V_6$ (qui est la même que celle de $-\alpha$).
     L'application $s \mapsto   \tilde{s}$ est injective. Comme le nombre d'éléments   de $V_6\sm\{0\}$ est 63, on obtient ainsi {\it une bijection de l'ensemble des
     réflexions de $G$ sur $V_6\sm\{0\}$}.  En particulier, l'ensemble des réflexions a {\it une structure naturelle de $\F_2$-espace projectif de dimension} 5.
     
     Deux réflexions  $s_1, s_2$ commutent si et seulement si $\tilde{s_1}$ et $\tilde{s_2}$ sont orthogonaux relativement à la forme alternée $a$ de $V_6$.
     
       Soit $C$ un cube maximal et soit $S$ sa base. Soit $\tilde{S}$ l'image de $S$
       dans $V_6$ par la bijection $s \mapsto  \tilde{s}$. Les 7 éléments de  $\tilde{S}$
       sont orthogonaux entre eux. Ils engendrent donc un sous-espace totalement
       isotrope de $V_6$ ayant au moins 7 éléments non nuls, donc de dimension $\geqslant 3$. Comme les sous-espaces totalement isotropes de $V_6$ sont de dimension $\leqslant 3$, on voit que {\it les bases des cubes maximaux sont les
       éléments non nuls d'un sous-espace totalement isotrope maximal de $V_6$}.
  En particulier:
  
    (a) {\it Les cubes maximaux correspondent bijectivement aux sous-espaces
    totalement isotropes maximaux de $V_6$}. Leur nombre est 135.
        
    (b) {\it La base d'un cube maximal a une structure naturelle de plan projectif sur le corps $\F_2$.}

    On a également:
    
    (c) {\it Le groupe $\Phi$ associé à un cube maximal est le groupe des permutations
   de sa base qui respectent sa structure de plan projectif sur $\F_2;$ c'est un groupe simple d'ordre $168$, isomorphe à} $\GL_3(\F_2)$. 
   
   \smallskip
   
   \noi {\it Démonstration de} (c). D'après le théorème de prolongement de Witt (cf. [Bo 59, cor.2 au th.1 du $\S$4), tout automorphisme
   linéaire d'un sous-espace totalement isotrope  de $V_6$ est la restriction
   d'un élément de $\Sp(V_6)$. D'où (c).
   
   \smallskip
   
     D'après la prop. 4.19, (c) permet d'interpréter les classes de conjugaison des involutions de $E_7$ comme les parties d'un plan projectif sur $\F_2$ modulo les automorphismes de ce plan. En degrés 0, 1, 2, il n'y a qu'une seule classe. 
     
     Par contre, en degré 3
   il y a deux classes: l'une d'elles correspond aux trois points d'une droite, l'autre aux trois sommets d'un triangle. On les distingue de la façon suivante: si $u$ est une involution de degré 3, écrivons $u$ comme produit $s_\alpha s_\beta s_\gamma$ de trois
   racines $\alpha, \beta, \gamma$ deux à deux orthogonales; alors $u$ est du type ``droite'' si et seulement
   si $\frac{\alpha+\beta+\gamma}{2}$ appartient au réseau $P$, autrement dit si ses produits scalaires avec les racines sont tous pairs.
   
   \smallskip
   
   \noi {\it Exemples}. Avec les notations de [Bo 68], $\S$VI.4.11, prenons pour $\alpha, \beta, \gamma$
   les racines $\alpha_2, \alpha_5,\alpha_7$; les produits scalaires de $\alpha + \beta + \gamma$ avec $\alpha_1,...,\alpha_7$ sont respectivement $0,2,0,-2,2,-2,2$; comme ils sont pairs, cela montre que $s_\alpha s_\beta s_\gamma$ est du type ``droite''. Si l'on prend $(\alpha, \beta, \gamma)= 
  (\alpha_1, \alpha_2,\alpha_5)$, le produit scalaire de $\alpha + \beta + \gamma$ avec $\alpha_3$ est égal à $-1$, donc impair; cela montre que $s_\alpha s_\beta s_\gamma$ est du type ``triangle''.
     
     \smallskip
       Par adjonction (cf. 3.11), le nombre de classes de conjugaison des involutions de degré 4, 5, 6, 7 est respectivement 2, 1, 1, 1.
   
      \medskip

     \noi 5.6. {\it Type $E_8$ - ordre $2^{14}\!\cdot\!3^5\!\cdot\!5^2\!\cdot\!7$} - [ATLAS], p. 85.
      
 \smallskip
 Le nombre des réflexions est $(248-8)/2=120$; les groupes $G_s^{^+}$ correspondants sont de type $ E_7$. Comme
 $E_7$ n'a qu'une seule classe d'involutions de degré $0, 1$ ou $2$, on en déduit la même chose pour $E_8$ jusqu'au degré 3, et, par adjonction, pour les degrés 5,6,7,8. Comme $E_7$ a deux classes d'involutions de degré 3 on voit que $E_8$ a au plus deux classes
 d'involutions de degré 4 (nous verrons ci-dessous qu'il en a deux). Comme $E_7$
 a 135 cubes maximaux, $E_8$ en a $120\times135/8 =  2025=3^4\!\cdot\!5^2;$ le groupe $\Phi$ associé à un tel cube est d'ordre  $2^{14}\!\cdot\!3^5\!\cdot\!5^2\!\cdot\!7/(2^8\!\cdot\!3^4\!\cdot\!5^2)= 2^6\!\cdot\!3\!\cdot\!7=8\!\cdot\!168$. Nous verrons
 dans (iii) ci-dessous quelle est sa structure.
 
 \smallskip
     Soit  $R$ le réseau des racines, et soit $V_8 = R/2R$. La fonction $\alpha \mapsto\frac{1}{2}$$\alpha \! \cdot \!\alpha$ définit par passage au quotient une forme quadratique non dégénérée $q$ sur le $\F_2$-espace vectoriel $V_8$. Cette forme est hyperbolique, autrement dit isomorphe
    à la forme $x_1x_2+x_3x_4+x_5x_6+x_7x_8$. Soit $Q$ l'ensemble des $x\in V_8$ tels que $q(x)=1$. On a $|Q| = 120$. Toute réflexion $s = s_\alpha$ définit un
élément  $\tilde{s}$ de $Q$, à savoir l'image de $\alpha$ (qui est aussi celle de $-\alpha$)
dans $V_8$. {\it L'application $s\mapsto \tilde{s}$ est une bijection de l'ensemble des réflexions sur
l'ensemble} $Q$. Deux réflexions distinctes $s$ et $s'$ commutent si et seulement si $\tilde{s}$ et $\tilde{s}'$ sont orthogonaux pour la forme $q$, autrement dit si la droite (géométrique) contenant $\tilde{s}$ et $\tilde{s}'$ est contenue dans la quadrique affine
d'équation $q(x)=1$.

\smallskip

  Soit $C$ un cube maximal, et soit $S$ sa base. Si $A\subset S$, 
 définissons $e(A) \in V_8$ par $e(A) = \sum_{s\in  A} \tilde{s}$. On a:
  
  \smallskip
  \noi (i) \ {\it Pour toute partie $A$ de $S$ à trois éléments, il existe
  un unique $s \in S\sm A$ tel que $\tilde{s}=e(A).$}
  
  \smallskip
  \noi {\it Démonstration}. Soient $x,y,z$ les éléments de $A$, et soit $a = e(A) = \tilde{x}+\tilde{y}+\tilde{z}$. Comme  $q(\tilde{x})=q(\tilde{y})=q(\tilde{z})=1$, et que $x,y,z$ sont deux à deux orthogonaux, on a $q(a)=3 = 1$, d'où $a\in Q$. Il existe donc une réflexion $s$ et une seule telle
  que $\tilde{s} = a$. Puisque $\tilde{s}$ est orthogonal aux $\tilde{z}$ pour $z\in S$, la réflexion  $s$ commute aux éléments de $S$, donc appartient à $S$, puisque  $S$  est la base d'un cube maximal.
 De plus, $s$ est distincte de $x$, car sinon, on aurait $\tilde{y}=\tilde{z}$, d'où  $y=z$; le même argument montre que
 $s$ est distincte de $y$ et de $z$. On a donc $s \in S\sm A$. 
  
\smallskip
Soit $\bb$ l'ensemble des parties $B$ de $S$ à quatre éléments telles que $e(B)=0$.
On peut reformuler (i) comme:

\smallskip 

\noi (ii) {\it Toute partie de $S$ à trois éléments est contenue dans un élément de~$\bb$ et un seul. }

\noi {\small (Autrement dit, le couple $(S,\bb)$ est un {\it système de Steiner} $S(3,4,8)$, cf.
[Wi 38].)}

\smallskip
Soit $\Phi$ le groupe de permutations de $S$ défini par le normalisateur de $C$, cf. 4.17. Son action sur
les parties à $4$ éléments a au moins deux orbites: une dans $\bb$ et une dans son complémentaire; comme il y a au plus deux orbites,
cela montre que $\Phi$ opère transitivement à la fois sur $\bb$ et sur son
complémentaire. En fait, on a un résultat plus précis:

\smallskip
\noi (iii) {\it Il y a sur $S$ une structure naturelle d'espace affine de dimension $3$ sur $\F_2$ . Le groupe $\Phi$ est le groupe d'automorphismes de cette structure$;$
il est isomorphe au produit semi-direct de $\GL_3(\F_2)$ par $\F_2^3$.}

\smallskip
\noi   {\it Démonstration}. Une partie $X$ d'un espace vectoriel sur $\F_2$ est un sous-espace affine si et seulement si elle a la propriété suivante:

 (Aff) $x,y,z \in X \ \Longrightarrow \ x+y+z \in X.$
 De plus, il suffit de vérifier (Aff) lorsque  $x,y,z$ sont distincts; en effet, si par exemple $x=y$, on a $x+y+z = z$, qui appartient à $X$.
   
     Appliquons ce critère à  l'ensemble $\tilde{S}$ des
  $\tilde{s}$ pour $s \in S$. La propriété (Aff) est satisfaite: en effet, si $\tilde{s_1}, \tilde{s_2}, \tilde{s_3}$ sont trois éléments distincts de $\tilde{S}$, d'après (i) il existe
  $s_4 \in S$ tel que $ \tilde{s_4} = \tilde{s_1}+ \tilde{s_2}+ \tilde{s_3}$. Il en résulte
  que $\tilde{S}$ est un sous-espace affine de $V_8$; d'où une structure naturelle
  d'espace affine sur $S$.
  
  \smallskip
  \noi(Cette structure affine résulte aussi
du théorème d'unicité des systèmes de Steiner $S(3,4,8)$, prouvé dans [Wi 38]; en effet,
un tel système est défini comme l'ensemble des sous-ensembles affines à 4 éléments d'un $\F_2$-espace affine de dimension 3.)

\smallskip
L'assertion
sur $\Phi$ en résulte. En effet, $\Phi$ est contenu dans le groupe d'automorphismes
de l'espace affine $S$, qui est produit semi-direct de $\GL_3(\F_2)$ par $\F_2^3$,
donc d'ordre $2^3\!\cdot\!168$, ce qui est aussi l'ordre de $\Phi$.
    
   \medskip
   
   D'après la prop. 4.19, on peut identifier les classes de conjugaison des involutions de $E_8$ aux parties d'un espace affine de dimension 3 sur $\F_2$, modulo
   automorphismes. En degré $\leqslant 3$, il n'y a qu'une seule classe. En degré 4, il y en a deux: celle (égale à $\bb$) des sous-espaces affines de dimension 2 (rectangles), et celle des parties à 4 éléments affinement
   indépendants (tétraèdres). On les distingue de la manière suivante : si $u$ est une
   involution de degré 4, on l'écrit comme produit de quatre réflexions $s_\alpha, s_\beta, s_\gamma, s_\delta$ commutant deux à deux; alors $u$ est de type
   ``rectangle'' si et seulement si $\alpha +\beta+\gamma+\delta$ appartient à $2R$.
   
   \smallskip
   
   \noi {\it Exemple}. Avec les notations de [Bo 68], $\S$VI.4.12, soit $\tilde{\alpha}$ la plus grande racine de $E_8$, et soit $\tilde{\beta}$ celle de son
   sous-système standard de type $E_7$. Prenons $(\alpha, \beta, \gamma, \delta) = (\tilde{\alpha}, \tilde{\beta}, \alpha_2, \alpha_3)$; l'involution $u = s_\alpha s_\beta s_\gamma s_\delta$ correspondante est de type ``rectangle". Par contre, si l'on remplace $(\alpha_2, \alpha_3)$ par $ (\alpha_5,\alpha_7)$, on obtient une involution de type ``tétraèdre''.

  \medskip 
     \noi 5.7. {\it Type $F_4$ - ordre $2^7\!\cdot\!3^2$.}
     
     \smallskip
        Le groupe $F_4$ contient $B_4$ comme sous-groupe d'indice $3$
     et $D_4$ comme sous-groupe normal d'indice 6. On a $F_4/D_4 \simeq \Sym_3$.
     L'action de $F_4$ sur $D_4$ par conjugaison définit un isomorphisme $F_4/\{\pm1\} \simeq \Aut(D_4)$. 
     
     Il y a deux classes de réflexions: 12 courtes
     et 12 longues. Elles sont échangées par l'automorphisme externe d'ordre 2:
     $ \{\alpha_1, \alpha_2,\alpha_3,\alpha_4\} \mapsto \{\alpha_4, \alpha_3,\alpha_2,\alpha_1\}$, où les $\alpha_i$ sont:
     
   \smallskip
  \hspace{25mm}$\alpha_{_1} \hspace{6mm}\alpha_{_2}\hspace{5mm}\alpha_{_3}\hspace{6mm} \alpha_{_4}$

\hspace{25mm}$\circ$-------$\circ$===$\circ$-------$\circ$.

\smallskip

 Les réflexions longues sont celles de $D_4$; les courtes s'en déduisent par l'automorphisme ci-dessus.

  Il y a deux classes de conjugaison d'involutions de degré 2. Pour l'une
  d'elles (par exemple $u =s_{\alpha_1}s_{\alpha_3}$) il existe un seul cube
  de rang 2 d'extrémité $u$; sa base est formée d'une réflexion courte et
  d'une réflexion longue. Une involution $u$ de l'autre classe est l'extrémité de deux cubes: un cube dont la base est formée de deux réflexions courtes, et un
  cube dont la base est formée de deux réflexions longues; le groupe $G_u^{^-}$ engendré
  par ces cubes est de type $B_2.$

  Comme pour $B_4$, il y a trois classes de conjugaison de cubes maximaux,
  caractérisées par le nombre de réflexions  courtes. Ce nombre
  est 4, 2 ou 0; le nombre des cubes correspondants est 3, 18 et 3.
  
\medskip

\noi 5.8. {\it Type $G_2$ - ordre $2^2\!\cdot\!3$.}

\smallskip

  C'est un groupe diédral d'ordre 12, isomorphe à $I_2(6)$. Il a 6
  réflexions, qui forment deux classes de conjugaison à 3 éléments chacune
  (correspondant aux racines courtes et aux racines longues).
  Il y a 3 cubes maximaux; leurs bases sont formées d'une réflexion courte et d'une réflexion longue; les groupes $\Phi$ correspondants sont triviaux.

\hspace{82mm} {\small 5}

\noi 5.9. {\it Type $H_3$ - ordre $2^3\!\cdot\!3\!\cdot\!5$.} Graphe de Coxeter: \ $\circ$-----$\circ$-----$\circ$

\smallskip

  On a $H_3 = \Alt_5 \times \{\pm1\}$, le plongement de $\Alt_5$ dans $\SO_3(\R)$ étant celui donné par l'action de $\Alt_5$ sur un icosaèdre régulier.

Il y a $3\!\cdot\!5=15$ réflexions, de la forme
  $(u,-1)$ où $u$ est un élément d'ordre 2 de $\Alt_5$; les groupes $G_s^{^+}$ correspondants sont de type $A_1 \times A_1$. Pour tout $d \leqslant 3$, il y a une  seule classe d'involutions de degré $d$. Il y a 5 cubes maximaux. Les groupes $\Phi$ correspondants sont d'ordre 3.
   
        \hspace{93mm} {\small 5}

\noi 5.10. {\it Type $H_4$ - ordre $2^6\!\cdot\!3^2\!\cdot\!5^2$.} Graphe de Coxeter: \ $\circ$-----$\circ$-----$\circ$----$\circ$

\smallskip

On peut définir $H_4$ de la manière suivante:
  Soit $\Alt_5^*$ le {\it groupe icosaédral binaire}, i.e. l'unique (à isomorphisme près)  extension non triviale de $\Alt_5$ par un groupe d'ordre 2. Soit $B = \Alt_5^*\times \Alt_5^*$. Soit $C = \{1,\sigma\}$ un groupe d'ordre 2, et soit $H =B\!\cdot\!C$ le produit semi-direct de $C$ et de $B$, où $C$ agit sur $B$ en permutant les deux facteurs.
  Soit $H'$ le quotient de $H$ par le sous-groupe central diagonal d'ordre 2 de $B$. Le groupe $H'$ a un plongement naturel dans le
  groupe orthogonal $O_4(\R)$, qui en fait un groupe de Coxeter irréductible
  de rang 4 de type $H_4$, cf. th.6.12. Il a 60 réflexions; ce sont les images dans $H'$
  des 120 éléments de la forme  $(a,a^{-1})\!\cdot\!\sigma$; les groupes $G_s^{^+}$ correspondants sont de type $H_3$.
  
 Pour tout $d \leqslant 4$, il y a une seule classe d'involutions de degré $d$.
Il y a 450 involutions de degré 2; d'après le th.4.5, chacune est l'extrémité d'un seul cube; leurs centralisateurs sont d'ordre 32. Il y a 60 involutions de degré 3; chacune est l'extrémité de 5 cubes. Il y a 75 cubes maximaux; les groupes $\Phi$ associés sont isomorphes à $\Alt_4$.

\medskip

\noi 5.11. {\it Type $I_2(m)$ - ordre $2m$}.
   
   \smallskip
   C'est un groupe diédral d'ordre $2m$. Il contient $m$ réflexions. C'est un groupe de
   Coxeter irréductible de rang 2 si $m>2$. On a $I_2(1) \simeq A_1$ et $I_2(2) \simeq A_1 \times A_1$.
   
   Si $m$ est pair, $I_2(m)$ contient $-1$, et il y a deux classes de réflexions. Le centralisateur d'une réflexion est un cube maximal. Si $m$ est divisible par $4$, le groupe $\Phi$ associé à un cube maximal est d'ordre 2 ; la base du cube est formée de deux réflexions conjuguées; il y a deux classes de cubes maximaux, qui ont chacune $m/4$ éléments; c'est le cas de $B_2 \simeq I_2(4)$. Si $m \equiv 2$ (mod 4), il y a $m/2$ cubes maximaux, qui sont conjugués entre eux; la base d'un tel cube est formée de deux réflexions non conjuguées; le groupe $\Phi$ est trivial; c'est le cas de $G_2\simeq I_2(6)$.
   
   Si $m$ est impair, le rang réduit est 1. Il y a une seule classe de réflexions; c'est le cas de $A_2\simeq I_2(3)$.
   
   \medskip 
   
   \noi 5.12.  {\it $\cc$-sous-groupes de petit indice.}
   
   \smallskip
   
     Dans les applications des cubes aux invariants cohomologiques ([Se 18]), il est commode de
     disposer de $\cc$-sous-groupes d'indice non divisibles par des nombres premiers comme 2, 3, 5.
    Voici
    les plus utiles:

    \medskip
    
    D'abord ceux d'indice impair:

    \noi 5.13. $ D_5 \subset E_6$ \quad \quad \quad: indice $3^3$; 
    
    \noi 5.14. $A_1 \times D_6 \subset E_7$: indice $3^2\!\cdot\!7$;
    
    \noi 5.15. $D_8 \subset E_8$  \quad \quad  \quad: indice $3^3\!\cdot\!5$;
    
    \noi 5.16. $B_4 \subset F_4$ \quad \quad  \ \: indice $3$;
    
    \noi 5.17. $A_1 \times A_1 \subset G_2$: indice $3$;
    
    \noi 5.18. $A_1 \times A_1 \times A_1 \subset H_3$: indice $3\!\cdot\!5$;
    
    \noi 5.19. $D_4 \subset H_4$ \quad \quad  \: indice $3\!\cdot\!5^2$. 
    
    \medskip
    
    Et aussi:
    
    \medskip
    
    \noi 5.20. $A_1 \times A_5 \subset E_6$: indice $2^2\!\cdot\!3^2;$

  \noi 5.21. $A_7 \subset E_7$ \quad \quad \quad: indice $2^3\!\cdot\!3^2;$
  
\noi 5.22. $E_6 \subset E_7$ \quad \quad \quad: indice $2^3\!\cdot\!7;$

   \noi 5.23. $A_4 \times A_4 \subset E_8$: indice $2^8\!\cdot\!3^3\!\cdot\!7;$
   
    \noi 5.24. $A_8 \subset E_8$ \quad \quad \quad: indice $2^7\!\cdot\!3\!\cdot\!5;$
    
    \noi 5.25. $A_2 \times E_6 \subset E_8$: indice $2^6\!\cdot\!5\!\cdot\!7;$
    
    \noi 5.26. $A_2 \times A_2 \subset F_4$: indice $2^5;$

\noi 5.27. $A_2 \subset G_2$ \quad \quad \quad: indice $2;$

\noi 5.28. $A_2 \subset H_3$ \quad \quad \quad: indice $2^2\!\cdot\!5;$

\noi 5.29. $I_2(5) \subset H_3$ \quad \quad: indice $2^2\!\cdot\!3;$

\noi 5.30. $A_2 \times A_2 \subset H_4$: indice $2^4 \!\cdot\!5^2  ;$

\noi 5.31. $I_2(5) \times I_2(5) \subset H_4$ \: indice $2^4\!\cdot\!3^2.$

\medskip

Pour les types crystallographiques, la plupart de ces inclusions  s'obtiennent en supprimant certains sommets du diagramme de Coxeter complété ([Bo 68], $\S$VI.4.3). Ainsi, pour le type $E_8$ ci-dessous,
on obtient 5.15, 5.23, 5.24 et 5.25 en supprimant respectivement $\alpha_1,
\alpha_5, \alpha_2$ et $\alpha_7$ dans le diagramme:

  \bigskip
  
\hspace{10mm}$\alpha_{_1} \hspace{6mm}\alpha_{_3}\hspace{5mm}\alpha_{_4}\hspace{6mm} \alpha_{_5}\hspace{5mm} \alpha_{_6}\hspace{5mm} \alpha_{_7}\hspace{5mm} \alpha_{_8}\hspace{6mm}- \!\tilde{\alpha}$

\hspace{10mm}$\circ$-------$\circ$------$\circ$-------$\circ$------$\circ$------$\circ$-------$\circ$ $^{.\ .\ .\ . }$ $\circ$

\hspace{28mm} $\mid$

\hspace{28mm} $\circ$  $\alpha_{_2}$

\bigskip

 \noi  Pour les inclusions 5.19, 5.30 et 5.31, relatives au type $H_4$, voir 6.14.
 
    \bigskip

 {\bf $\S$6. Groupes de Coxeter de rang $4$}. 
 
 \smallskip
 
 Les sous-groupes finis de $\SO_3(\R)$ peuvent être utilisés pour construire des groupes de Coxeter de rang 4 par un procédé de
 ``dédoublement''. Avant de définir cette construction (ce qui sera fait dans 6.11),
  un certain nombre de préliminaires sont nécessaires.
  
  \medskip\medskip

\noi 6.1. {\it  Premier dédoublement $:$ le groupe $\Gamma_2$.}
  
    \medskip
    
    Soit $\Gamma$ un groupe. Soit $\Gamma_2$ le sous-groupe
    de $\Gamma^2=\Gamma \times \Gamma$ engendré par les $(x,y)$ tels que $xy=1$.
 
 \medskip
 
 \noi {\bf Proposition 6.2}. {\it Un élément $(a,b)$ de $\Gamma^2$ appartient à $\Gamma_2$ si et seulement si $ab \equiv 1 $ \ {\rm mod} $D(\Gamma)$.}
 
 \medskip
 \noi {\bf Corollaire 6.3}. {\it Le groupe $\Gamma_2$ est un sous-groupe normal
 de $\Gamma^2$. L'application} $(x,y) \mapsto x$ mod $D(\Gamma)$ {\it définit un isomorphisme  $\Gamma^2/\Gamma_2$ sur  $\Gamma^{\rm ab}$.}
 
 \smallskip
 \noi En particulier, on a $\Gamma_2=\Gamma^2$ lorsque $\Gamma^{\rm ab}=1.$
 
 \medskip
 
 \noi {\it Démonstration de la prop.} 6.2. La relation:
 
 \medskip
 
 \ \   $(xyx^{-1}y^{-1},1) = (x,x^{-1})\cdot(yx^{-1},xy^{-1})\cdot(y^{-1},y)$
 
 \medskip
 
\noi entraîne que $\Gamma_2$ contient $D(\Gamma) \times 1$, donc aussi $D(\Gamma)\times D(\Gamma)$; cela montre en particulier que $\Gamma_2$ est un groupe normal de $\Gamma^2$. De plus, l'image de $\Gamma_2$ dans $\Gamma^{\rm ab} \times \Gamma^{\rm ab}$ est l'antidiagonale de ce groupe, i.e.
 l'ensemble des couples $(x,x^{-1})$. D'où 6.2 et 6.3.
 
 \medskip
 \noi 6.4. {\it Deuxième dédoublement $:$ les groupes } $\Gamma^2_\sigma$ et $\Gamma_{2,\sigma}$.
 
 \medskip
   Soit $\{1,\sigma\}$ un groupe d'ordre 2, que l'on fait opérer sur $\Gamma^2$
   par $(x,y)\mapsto (y,x)$. Nous noterons $\Gamma^2_\sigma$ (resp. $\Gamma_{2,\sigma}$) le produit semi-direct de $\{1,\sigma\}$ par $\Gamma^2$
 (resp. par  $\Gamma_2$)\footnote{Autrement dit, $\Gamma^2_\sigma$ est le {\it produit en couronne} $\Gamma \wr \{1,\sigma\}$ de $\Gamma$ par $ \{1,\sigma\}$.}. Le groupe $\Gamma_{2,\sigma}$ est un sous-groupe normal de $\Gamma^2_\sigma$, et le quotient $\Gamma^2_\sigma/\Gamma_{2,\sigma}$ est isomorphe à $\Gamma^{\rm ab}$.
 
 \smallskip
   Soit $a \in \Gamma$, et soit $s_a = (a,a^{-1})\sigma \in \Gamma_{2,\sigma}$. On a:  $$s_a^2 = (a,a^{-1})\sigma (a,a^{-1})\sigma = (a,a^{-1})(a^{-1},a) \sigma \sigma= 1.$$
\noi Ainsi, $s_a$ est un élément d'ordre 2 de $\Gamma^2_\sigma$. On a $s_1=\sigma$. 
 \medskip
 
 \noi {\bf Proposition 6.5}. {\it Le groupe $\Gamma_{2,\sigma}$  est engendré par les $s_a, a \in \Gamma$.}
 
 \smallskip
 
  En effet, on a $s_as_1= (a,a^{-1})$, et les éléments $(a,a^{-1})$ engendrent $\Gamma_2$.
  
  \medskip
    
\noi 6.6.   {\it Dernière étape$:$ les groupes $B(\Gamma)$ et $B(\Gamma)^c$.}
  
  \medskip
A partir de maintenant, {\it on se donne un élément central $e$ d'ordre $2$ de} $\Gamma$; on pose $C = \{1,e\}$ et l'on note $\Gamma_0$ le quotient $\Gamma/C$.     
  
  Soit $C_2 = \{1,(e,e)\}$ le sous-groupe diagonal de $C\times C$. C'est un sous-groupe central de $\Gamma_{2,\sigma}$. Le quotient de $\Gamma^2_\sigma$ par $C_2$ sera noté $B(\Gamma)$, et le quotient de $\Gamma_{2,\sigma}$
  par $C_2$ sera noté $B(\Gamma)^c$.
  
  \smallskip
  
  Ces groupes  sont  produits semi-directs de $\{1,\sigma\}$ par $\Gamma^2/C_2$ et $\Gamma_2/C_2$ respectivement. 
  
  \smallskip
  
  Lorsque $\Gamma$ est fini, on a $| B(\Gamma)| = | \Gamma |^2$
   et $| B(\Gamma)^c| = |\Gamma|^2/| \Gamma^{\rm ab}|$.

\medskip
 \noi 6.7. Pour tout $a\in \Gamma$, soit $\sigma_a$ l'image de $s_a$
  dans $B(\Gamma)$; on a $\sigma_a \in B(\Gamma)^c,$ et, d'après la prop.6.5, {\it les
  $\sigma_a$ engendrent} $B(\Gamma)^c$. On a $\sigma_{ae} = \sigma_a$,
  puisque $(a,a^{-1})$ et $(ae,a^{-1}e)$ sont congrus mod $C_2$; on peut donc
  considérer que les $\sigma_a$ ont un sens pour $a \in \Gamma_0$.
  
  \medskip
  
  \noi {\bf Proposition 6.8}. {\it Soient $a,b \in \Gamma_0$. L'ordre de $\sigma_a\sigma_b$ est 
  égal à l'ordre de $ab^{-1}$.}
  
  \smallskip 
  
  Soient $a'$ et $b'$ des représentants de $a$ et $b$ dans $\Gamma$. L'ordre de $ab^{-1}$ est le plus entier $n > 0$ tel que $(ab^{-1})^n=1$ ou $e$. Le produit
  $\sigma_a\sigma_b$ est l'image de $s_{a'}s_{b'}$ dans $\Gamma^2/C_2$. Or on a:
  $$s_{a'}s_{b'}=  (a',a'^{-1})\sigma(b',b'^{-1})\sigma = (a',a'^{-1})(b'^{-1},b')= (a'b'^{-1},a'^{-1}b').$$ 
  D'où $(s_{a'}s_{b'})^m = ((a'b'^{-1})^m,(a'^{-1}b')^m)$ pour tout $m$. Comme $(a'b'^{-1})^m$ et $(a'^{-1}b')^m$ sont conjugués, si l'un d'eux est égal à $1$ (resp. à ${e}$), l'autre l'est aussi. Le plus petit $m$ pour lequel cela ait lieu est donc $n$. Cela montre que l'ordre de $s_{a'}s_{b'}$ modulo $C_2$ est égal à $n$; autrement dit, l'ordre de $\sigma_a\sigma_b$ est $n$.

 \bigskip
 \noi 6.9. {\it Exemple de groupe $\Gamma:$ le groupe des quaternions de norme $1$.}
 
 \bigskip
 
   Soit $\H$ le corps des quaternions sur $\R$, et soit $\s_3$ le groupe des éléments $z$ de $\H^\times$ tels que $z\overline{z}=1$.  On a $\s_3 \simeq \SU_2(\C)$.  Choisissons pour $\Gamma$ le groupe $\s_3$, et pour $e$ l'élément -1; on a $C = \{1,-1\}$ et $\Gamma_0 \simeq \SO_3(\R)$. 
   
   \smallskip
   
  \noi {\bf Proposition 6.10}. {\it Le groupe $\Gamma^2/C_2 = \Gamma_2/C_2$ est isomorphe à $\SO_4(\R)$ et le groupe  $B(\s_3)^c$ est isomorphe à} $\O_4(\R)$.

  \smallskip
 \noi {\it Démonstration}.  On obtient un isomorphisme $\varphi: \Gamma^2/C_2 \to \SO_4(\R)$ en associant à un 
   élément $(a,b)$ de $\Gamma^2$ l'application de $\H \simeq \R^4$ dans lui-même donnée par $z \mapsto az\overline{b}$. On prolonge $\varphi$ à $B(\Gamma)^c= B(\s_3)^c$ en définissant $\varphi(\sigma)$  comme la réflexion 
$z \mapsto - \overline{z}$; on obtient ainsi  {\it un isomorphisme} $ \varphi:B(\s_3)^c \simeq \O_4(\R)$.
    
    \smallskip
   Cet isomorphisme  transforme les $\sigma_a$ en les réflexions de $\O_4(\R)$:
 plus précisément,  $\sigma_a: z \mapsto -a\overline{z}a$ est l'unique réflexion qui
 transforme $a$ en -$a$. 
 
 \medskip
 
 \noi 6.11. {\it Exemples de groupes $\Gamma:$ les groupes de Coxeter associés aux groupes finis de rotations de $\R^3$.}
 
 \medskip
 
   Soit $\Gamma_0$ un sous-groupe fini de $\SO_3(\R)$, et soit $\Gamma $ son image réciproque dans $\s_3$ (``groupe double"),
   muni de l'élement central $e=-1$ comme ci-dessus. D'après ce qui précède,
   le groupe $B(\Gamma)^c$ est un sous-groupe de Coxeter de $B(\s_3)^c = \O_4(\R)$. De façon plus précise:

   \medskip
   
   \noi {\bf Théorème 6.12}. 
   
   (a) {\it Si $\Gamma_0$ est cyclique d'ordre $m$, $B(\Gamma)^c$ est de type} $I_2(m)$.
  
   (b)  {\it Si $\Gamma_0$ est diédral d'ordre $2m$, $B(\Gamma)^c$ est de type} $I_2(m)\times I_2(m)$.
   
   (c)  {\it Si $\Gamma_0$ est isomorphe à $\Alt_4$, \ \ $B(\Gamma)^c$ est  de type} $D_4$.
   
 (d)  {\it Si $\Gamma_0$ est isomorphe à $\Sym_4$, $B(\Gamma)^c$ est  de type} $F_4$.

 (e)  {\it Si $\Gamma_0$ est isomorphe à $\Alt_5$, \ \ $B(\Gamma)^c$ est  de type} $H_4$.

 \medskip
 
 \noi {\it Démonstration}. Dans les cas (a),...,(e), on a  $|\Gamma| = 2m, 4m, 24, 48, 120$, et $|\Gamma^{\rm ab}| = 2m, 4, 3, 2, 1$, d'où $| B(\Gamma)^c| = 2m, 4m^2, 2^6\!\cdot\!3, 2^7\!\cdot\!3^2, 2^6\!\cdot\!3^2\!\cdot\!5^2$, qui sont
 les ordres de $I_2(m),..., H_4$. De plus, d'après 6.8, l'ordre maximum d'un produit
 de deux réflexions est respectivement $m, m, 3, 4, 5$. Ces renseignements
 suffisent à prouver (a), ..., (e). [Noter que le cas (e) est traité dans  [Bo 68], $\S$VI.4, exerc.12.]
 
 \medskip
 
 \noi 6.13. {\it Bases explicites}.

\medskip

 On peut préciser (et redémontrer) le th. 6.12 en donnant une présentation de Coxeter pour chacun des  groupes $B(\Gamma)^c$. Voici des choix possibles pour chaque type (a), ..., (e):
 
 \medskip
 {\it Type} (a). Le groupe $\Gamma_0$ est cyclique d'ordre $m$. Soit $x$ un générateur. D'après 6.8, l'ordre de $\sigma_1\sigma_x$ est $m$. Les
 deux réflexions $\sigma_1$ et $\sigma_x$ engendrent un sous-groupe diédral
 de $B(\Gamma)^c$ d'ordre $2m$. Comme $| B(\Gamma)^c| = 2m$ d'après 6.6,
 on en conclut que $B(\Gamma)^c$ est isomorphe à $I_2(m)$.
 
 \medskip
 
 {\it Type} (b). Dans ce cas, $\Gamma_0$ est engendré par deux éléments $s$ et $x$, où $s$ est d'ordre 2, $x$ est d'ordre $m$ et $sxs = x^{-1}$. Le groupe $B(\Gamma)^c$ contient les
 quatre réflexions $\sigma_1, \sigma_x, \sigma_s, \sigma_{sx}$. D'après 6.8,
 l'ordre du produit $\sigma_1\sigma_{sx}$ est égal à 2, ce qui signifie que $\sigma_1$ et $\sigma_{sx}$ commutent. Le même argument montre que $\sigma_1$ commute à
 $\sigma_s$, que $\sigma_x$ commute à $\sigma_s$ et $\sigma_{sx}$, et que
 les produits
 $\sigma_1\sigma_x$ et  $\sigma_s\sigma_{sx}$ sont d'ordre $m$. On en conclut
 que $B(\Gamma)^c$ contient le groupe de Coxeter de diagramme:

  \smallskip
\hspace{28mm} $_m$ \hspace{14mm} $_m$

\hspace{25mm}$\circ$-------$\circ$ \hspace{6mm}$\circ$-------$\circ$  \ \ .
  
\hspace{25mm}$\sigma_{_1} \hspace{6mm}\sigma_{_x}\hspace{6mm}\sigma_{_s}\hspace{7mm}\sigma_{_{sx}} $

\medskip

\noi Comme   l'ordre de $B(\Gamma)^c$ est $4m^2$, ce groupe est de type $I_2(m)\times I_2(m)$.
  
 \medskip
 {\it Type} (c). On a $\Gamma_0=\Alt_4$. Pour prouver que $B(\Gamma)^c$ est de type $D_4$, il suffit de trouver quatre
 réflexions $\sigma_a, \sigma_b, \sigma_c, \sigma_d$, avec $a,b,c,d \in \Alt_4$,  telles que chacun des
 produits de $\sigma_a$ avec les deux autres soit d'ordre 3 et que les produits de deux quelconques des trois autres soient d'ordre 2. D'après 6.8, cela revient à demander que $a^{-1}b, a^{-1}c, a^{-1}d$ soient d'ordre 3 et $b^{-1}c, c^{-1}d, d^{-1}b$ d'ordre 2. Les permutations suivantes répondent à ces conditions:   $  a= 1, b = (123), c = (142), d = (134)$; on a en effet $b^{-1}c = (23)(14), c^{-1}d= (13)(24), d^{-1}b= (12)(34)$.
 
 Noter que, si  $f= (243)$, les quatre réflexions $\sigma_b, \sigma_c,\sigma_d,\sigma_f$ sont la base d'un cube maximal.
 
\medskip
  
   {\it Type} (d). On a $\Gamma_0=\Sym_4$. On prend $a=1, b = (123), c= (14), d= (12)$, et l'on obtient le diagramme de $F_4$:
  
      \smallskip
\hspace{39mm} $_{4}$

\hspace{25mm}$\circ$-------$\circ$--------$\circ$-------$\circ$ .
  
\hspace{25mm}$\sigma_{_a} \hspace{6mm}\sigma_{_b}\hspace{6mm}\sigma_{_c}\hspace{7mm}\sigma{_d} $

\medskip

Si $f=(12)(34), g=(13)(24), h=(14)(23)$, les quatre réflexions $\sigma_1, \sigma_f,\sigma_g,\sigma_h$ sont la base d'un cube maximal.

  {\it Type} (e). On a $\Gamma_0=\Alt_5$. On prend $a=1, b= (12345), c = (15)(34), d= (15)(24)$ et l'on obtient
 le diagramme de  $H_4$:

  \smallskip
\hspace{29mm} $_5$

\hspace{25mm}$\circ$-------$\circ$--------$\circ$-------$\circ$ .
  
\hspace{25mm}$\sigma_{_a} \hspace{6mm}\sigma_{_b}\hspace{6mm}\sigma_{_c}\hspace{7mm}\sigma{_d} $

\bigskip

 \noi 6.14.     {\it Application à la construction de sous-groupes.} 
      
   \smallskip
       Le fait que $\Alt_4$ soit un sous-groupe de $\Alt_5$ entraîne
      une inclusion analogue pour les groupes $B(\Gamma)^c$ ; cela prouve
      que $D_4$ est isomorphe à un $\cc$-sous-groupe de $H_4$. Cette inclusion
     (ainsi que la méthode pour l'obtenir) m'a été signalée par C. Bonnafé; on la trouve aussi dans [KKAS 06] et [DPR 13].
     
       De même, le fait que $\Alt_5$ contienne des sous-groupes diédraux d'ordre 6
       et 10 montre que $H_4$ contient des $\cc$-sous-groupes de type $A_2 \times A_2$ et $I_2(5) \times I_2(5)$. Cela donne les inclusions 5.30 et 5.31.
       
       \newpage
    
           \begin{center}
      
        {\sc Références}
        
        \end{center}

\noi [ATLAS] J.H. Conway, R.T. Curtis, S.P. Norton, R.A. Parker \& R.A. Wilson, {\it Atlas of Finite Groups}, Clarendon Press, Oxford, 1985; second corrected edition, 2003.

\noi [Bo 59]  N. Bourbaki, {\it Algèbre, Chap. IX}, Hermann, Paris, 1959; English translation, Springer-Verlag, 1981.

\noi [Bo 68]  N. Bourbaki, {\it Groupes et alg\`ebres de Lie, Chap. IV-VI}, Hermann, Paris, 1968; English translation, Springer-Verlag, 2002.

\noi  [BS 20]   E. Bayer-Fluckiger \& J-P. Serre,  {\it Lines on cubic surfaces, Witt invariants and Stiefel-Whitney classes}, arXiv:1909.05312; à paraître dans Indag. Math.

\noi [Ca 72] R.W. Carter, {\it Conjugacy classes in the Weyl group}, Compos. Math. 25 (1972), 1-59.

\noi [CE 56] H. Cartan  \& S. Eilenberg, {\it Homological Algebra}, Princeton Univ. Press, 1956.

\noi [DPR 13] J.M. Douglass, G. Pfeiffer,  G. Röhrle, {\it On reflection subgroups of finite Coxeter groups},  Comm. Algebra 41 (2013), 2574–2592.

\noi [GH 21] S. Gille \& C. Hirsch, {\it On the splitting principle for cohomological invariants of reflection groups}, arXiv:1908.08146; à paraître dans Transformation Groups. 
 
\noi [Hi 10] C. Hirsch, {\it Cohomological invariants of reflection groups}, Diplomarbeit, L.M.U., München, 2009, arXiv: 1805.04670v1.
 
 \noi [Hi 20] C. Hirsch, {\it On the decomposability of mod $2$ cohomological invariants of Weyl groups}, Comm. Math. Helv. 95 (2020), 765-809.

\noi [Ho 88]  R.B. Howlett, {\it On the Schur multipliers of Coxeter groups}, J. London Math. Soc. 38 (1988), 263-276.

\noi [IY 65] S. Ihara \& T. Yokonuma, {\it On the second cohomology groups $($Schur multipliers$)$ of the finite Coxeter groups}, J. Fac. Sci. Univ. Tokyo 11 (1965), 155-171.

\noi [Ka 01] R. Kane, {\it Reflection Groups and Invariant Theory}, Springer-Verlag, 2001.

\noi [KKAS 06] M. Koca, R. Koç, M. Al-Barwani \& S. Al-Farsi, {\it Maximal subrgoups of the Coxeter group $W(H_4)$ and quaternions}, Linear Alg. and Appl. 412 (2006), 441-452.

\noi [La 96] S. Lang, {\it Topics in Cohomology of Groups}, L.N.M. 1625, Springer-Verlag, 1996.

\noi [Ma 74] Y. Manin, {\it Cubic Forms, Algebra, Geometry, Arithmetic}, North-Holland, 1974; second edit., 1986.

\noi [Sc 11] I. Schur, {\it Über die Darstellung der symmetrischen und der alternierenden Gruppen durch gebrochene lineare Substitutionen}, J. Crelle 139
 (1911), 155-250; {\it Ges. Abh.} I, 346-441.

\noi [Se 03] J-P. Serre, {\it Cohomological invariants, Witt invariants and trace forms}, notes by Skip Garibaldi, ULS 28,  AMS (2003), 1-100.

\noi [Se 18] J-P. Serre, {\it Cohomological invariants {\rm mod} $2$ of Weyl groups}, Oberwolfach reports 21 (2018), 1284-1286; arXiv:1805.07172.

\noi [Sp 74] \ T.A. Springer, {\it Regular elements of finite reflection groups}, Invent. math. 25 (1974), 159-198.

\noi [Sp 82] \  T.A. Springer, {\it Some remarks on involutions in Coxeter groups}, Comm. Algebra 10 (1982), 631-636.

\noi [Wi 38] E. Witt, {\it Über Steinersche Systeme}, Hamburg Abh.12 (1938), 265-275; {\it Ges. Abh.}, 288-298.

\vspace{6mm}

\noi Collège de France, Paris

\noi jpserre691@gmail.com     
\end{document}